\documentclass[microtype]{gtpart}   

\usepackage{color,graphicx,psfrag, verbatim,amssymb,amscd,enumerate,subfigure,ytableau,tikz,texdraw,xypic,appendix, accents}
\usepackage[all]{xy}

\title{The Gerby Gopakumar-Mari\~no-Vafa Formula}

\author{Dustin Ross} 
\givenname{Dustin}
\surname{Ross}
\address{Department of Mathematics\\Colorado State University\\\newline
         Fort Collins, CO 80523\\USA}
\email{ross@math.colostate.edu}

\author{Zhengyu Zong} 
\givenname{Zhenyu}
\surname{Zong}
\address{Department of Mathematics\\Columbia University\\\newline
         New York, NY 10027\\USA}
\email{zz2197@math.columbia.edu}

\subject{primary}{msc2010}{14N35}
\subject{primary}{msc2010}{53D45}
\subject{secondary}{msc2010}{05E05}

\arxivreference{1208.4342}  
\arxivpassword{kvr3m}

\newcommand{\re}{\mathrm{e}}
\newcommand{\cha}{\raisebox{.05cm}{$\chi$}}

\newcommand{\Tot}{\text{Tot}}
\newcommand{\Pic}{\text{Pic}}
\newcommand{\ord}{\text{ord}}

\newcommand{\proj}{\mathbb{P}}
\newcommand{\N}{\mathbb{N}}
\newcommand{\aut}{\mathrm{Aut}}

\newcommand{\M}{\overline{\mathcal{M}}}
\newcommand{\so}{\mathcal{O}}
\newcommand{\X}{\mathcal{X}}

\newcommand{\B}{\mathcal{B}}
\newcommand{\La}{\mathcal{L}}

\newcommand{\G}{\mathcal{G}}

\DeclareMathOperator{\Gcd}{gcd}

\newcounter{counter}

\newtheorem{dummy}{}[section]
\newtheorem{lemma}[dummy]{Lemma}

\newtheorem{theorem}[dummy]{Theorem}
\newtheorem{corollary}[dummy]{Corollary}
\newtheorem{reduction}[dummy]{Reduction}
\newtheorem{theorem1}[counter]{Theorem}
\newtheorem{theorem2}[counter]{Theorem}

\newtheorem{definition}[dummy]{Definition}

\newtheorem{remark}[dummy]{Remark}

\numberwithin{equation}{section}

\begin{document}

\begin{abstract}

We prove a formula for certain cubic $\Z_n$-Hodge integrals in terms of loop Schur functions.  We use this identity to prove the Gromov-Witten/Donaldson-Thomas correspondence for local $\Z_n$-gerbes over $\proj^1$.

\end{abstract}

\maketitle

\section{Introduction}

\subsection{Statement of Results}

The Gopakumar-Mari\~no-Vafa formula, proven independently by Liu, Liu, and Zhou \cite{llz:pmvf} and Okounkov and Pandharipande \cite{op:hiaiotu}, evaluates certain generating functions of cubic Hodge integrals on moduli spaces of curves in terms of Schur functions, a special basis of the ring of symmetric functions.  The formula can be interpreted as one instance of the Gromov-Witten/Donaldson-Thomas correspondence for Calabi-Yau (CY) 3-folds.  In this paper, we generalize the Gopakumar-Mari\~no-Vafa formula to $\Z_n$-Hodge integrals and we show that this formula can be viewed as one instance of the orbifold GW/DT correspondence.

In particular, we define generating functions $\tilde{V}_\mu^\bullet(a)$ of cubic $\Z_n$-Hodge integrals on moduli spaces of stable maps to the classifying space $\B\Z_n$.  These generating functions are indexed by conjugacy classes $\mu$ of the generalized symmetric group $\Z_n\wr S_d$ and are closely related to the GW orbifold vertex developed by the first author \cite{r:lgoa}.  In place of the Schur functions in the usual Gopakumar-Mari\~no-Vafa formula, we introduce generating functions $\tilde{P}_\lambda(a)$ which are specializations of \textit{loop} Schur functions, developed by Lam and Pylyavskyy \cite{lp:tpilg1} and further investigated by the first author \cite{r:cmn}.  These generating functions are indexed by irreducible representations $\lambda$ of $\Z_n\wr S_d$ and are closely related to the DT orbifold vertex developed by Bryan, Cadman, and Young \cite{bcy:otv}.  The main result is the following correspondence via the character values $\cha_\lambda(\mu)$ of $\Z_n\wr S_d$.

\begin{theorem1}\label{thm1}
After an explicit change of variables,
\[
\tilde{V}_\mu^\bullet(a)=\sum_\lambda \tilde{P}_\lambda(a)\frac{\cha_\lambda(\mu)}{z_\mu}
\]
\end{theorem1}

There are $n$ distinct $\Z_n$-gerbes $\G_k$ ($0\leq k<n$) over $\proj^1$ classified by $H^2(\proj^1,\Z_n)$.  We define $\X$ to be a \textit{local $\Z_n$-gerbe over $\proj^1$} if $\X$ is isomorphic to the total space of a rank two Calabi-Yau orbifold bundle over some $\G_k$.  Applying the gluing rules of \cite{r:lgoa} and \cite{bcy:otv}, Theorem \ref{thm1} leads to a proof of the orbifold GW/DT correspondence for local $\Z_n$-gerbes over $\proj^1$.

\begin{theorem2}\label{thm2}
After an explicit change of variables, the GW potential of any local $\Z_n$-gerbe over $\proj^1$ is equal to the reduced, multi-regular DT potential.
\end{theorem2}

This is the first example of the GW/DT correspondence for orbifold targets with nontrivial curve classes contained in the singular locus.

\subsection{Context and Motivation}

Atiyah-Bott localization \cite{ab:tmmaec} has proven to be an extremely powerful tool in both GW and DT theory of toric CY 3-folds.  In particular, it has led to the development of the (orbifold) \textit{topological vertex} \cite{akmv:tv,orv:qcycc,lllz:mttv,bcy:otv,r:lgoa}, a basic building block for the GW or DT theory of all toric CY 3-folds.  In the GW case the vertex can be computed as a generating function of (abelian) Hodge integrals, whereas in the DT case the vertex can be computed as a generating function of (colored) 3d partitions.

The topological vertex formalism provides us with an algorithm for proving conjectural correspondences related to GW and DT theory: first prove that the correspondence holds locally for the vertex, then show that it is consistent with the gluing laws.  In the smooth case, this approach was utilized to prove the GW/DT correspondence for toric 3-folds, initiated by Li, Liu, Liu and Zhou \cite{llz:pmvf,llz:ftphi,lllz:mttv} and concluding with the work of Maulik, Oblomkov, Okounkov, and Pandharipande \cite{moop:gwdtc}.

In orbifold Gromov-Witten theory, the first example of this local-to-global approach appeared in work of the first author and Cavalieri \cite{cr:ogwcrc} where it was used to prove an example of the crepant resolution conjecture.  It was further developed in \cite{r:lgoa}, where a correspondence between the $A_{n-1}$ GW and DT vertex theories was suggested.  In \cite{z:gmvf}, the second author proved this correspondence for the effective one-leg $A_{n-1}$ vertex.  One consequence of the results in \cite{z:gmvf} is the orbifold GW/DT correspondence for local footballs (orbifolds with coarse space $\proj^1$ and smooth away from $0$ and $\infty$).  The main focus of this paper is the \textit{ineffective} one-leg $A_{n-1}$ vertex.

In the ineffective case, several new challenges arise.  On the GW side one can no longer utilize the $\Z_n$-Mumford relation which was the key tool in \cite{z:gmvf}.  Moreover, the orbifold structure at the nodes of the source curve is no longer determined by the degree of the corresponding map, barring us from using the standard tools in the representation theory of the symmetric group.  In order to overcome the first challenge, instead of evaluating the abelian Hodge integrals explicitly, we develop a deterministic set of bilinear relations by localizing relative maps into cyclic gerbes over $\proj^1$.  The difficulties lie in developing an efficient set of useful relations and then showing that these relations are invertible.  As suggested in \cite{r:lgoa}, the latter obstacle is overcome by encoding the twisted partitions as conjugacy classes in the generalized symmetric group $\Z_n\wr S_d$.  Many of the combinatorial tools from the study of the representation theory of $S_d$ can then be generalized to $\Z_n\wr S_d$, and these tools are crucial in proving our main results.

On the DT side, the effective case can be interpreted in terms of Schur functions, but we lose this interpretation when we pass to the ineffective case.  However, we observe that the DT vertex can naturally be interpreted as specializations of \textit{loop} Schur functions.  The combinatorial structure of the loop Schur functions provides us in turn with useful properties of the DT vertex which are pivotal in the arguments of this paper.

Many interesting questions arise from this work.  First, the results of this paper give the first example of the orbifold GW/DT correspondence for a target which contains nontrivial curve classes which lie entirely in the singular locus.  In this case, it is necessary to discard a significant amount of information on the DT side by restricting to the multi-regular contributions.  It would be interesting to generalize orbifold GW theory to account for this extra data and one possible approach seems to lie in the \textit{very twisted stable maps} developed by Chen, Marcus, and \'{U}lfarsson \cite{cmu:vtsm}.  Secondly, since the current work completes the one-leg $A_{n-1}$ GW/DT correspondence, another natural extension of this work is to extend the results herein to the two, and ultimately the three-leg $A_{n-1}$ vertex.  Finally, the $A_{n-1}$ vertex is by far the easiest geometry in both GW and DT theory.  It would be extremely interesting to study if/how the GW/DT vertex correspondence extends to noncyclic and/or non hard-Lefschetz orbifolds.

\subsection{Acknowledgments}

We are grateful to Jim Bryan and Paul Johnson for enlightening conversations, insights, and suggestions along the way.  We are also pleased to acknowledge the referee for many helpful suggestions which greatly improved the exposition and in particular for suggesting that we merge our individual works into the current joint paper.  Most of all, we are greatly indebted to our advisors Renzo Cavalieri and Chiu-Chu Melissa Liu, respectively.  Their expertise and guidance has always been extremely enlightening and this project could not have been possible without them.

\subsection{Plan of the Paper}

After setting up notation and giving a precise statement of Theorems 1 and 2 in Section \ref{sec:background}, we study the geometry of the framed GW vertex $\tilde{V}_\mu^\bullet(a)$ in Section \ref{sec:geometry}.  In particular, we develop a set of bilinear equations relating the GW vertex to generating functions of certain rubber integrals.  In Section \ref{sec:wreath}, we interpret these rubber integrals in terms of wreath Hurwitz numbers and apply the Burnside formula to write the bilinear relations in terms of the characters of the generalized symmetric group $\Z_n\wr S_d$.  We then show in Section \ref{sec:invertible} and Section \ref{sec:linearalgebra} that these relations uniquely determine the GW vertex.  Sections \ref{sec:combinatorics}, \ref{sec:reptheory}, and \ref{sec:proof} are devoted to proving that the DT vertex also satisfies these bilinear relations.  In Section \ref{sec:combinatorics}, we recall the definition of loop Schur functions and the main results from \cite{r:cmn}.  We also recall a hook-length formula from \cite{er:cpgisf} and \cite{n:hfgt} which relates the loop Schur functions to the framed DT vertex $\tilde{P}_\lambda(a)$.  In Section \ref{sec:reptheory}, we study the representation theory of $\Z_n\wr S_d$ where the main tool is the wreath Fock space.  Finally, in Section \ref{sec:proof} we put everything together to prove Theorem \ref{thm1}.  In Section \ref{sec:app} we use gluing rules developed in \cite{r:lgoa} and \cite{bcy:otv} to show how the GW/DT correspondence for local $\Z_n$-gerbes over $\proj^1$ follows from Theorem \ref{thm1}.

\section{Background and Notation}\label{sec:background}

In this section we set up notation which will be used throughout the paper and we give a precise statement of the main results.

\subsection{Partitions}\label{sec:parts}

For each positive integer $n$ we fix a generator of the cyclic group \[\Z_n=\left\langle\xi_n:=e^{\frac{2\pi\sqrt{-1}}{n}}\right\rangle.\]  When no confusion arises, we write the generator simply as $\xi$.  It is well known that $n$-tuples of partitions naturally correspond to conjugacy classes and irreducible representations of $\Z_n\wr S_d$, see e.g. \cite{m:sfhp}.  We will use $\mu$ and $\nu$ to denote $n$-tuples of partitions corresponding to conjugacy classes and reserve $\lambda$ and $\sigma$ to refer to irreducible representations.  We let $\cha_\lambda(\mu)$ denote the value of the character of the irreducible representation $\lambda$ on the conjugacy class $\mu$.

Consider the $n$-tuple of partitions
\[
\mu=\left((d_1^0,...,d_{l_0}^0),...,(d_{1}^{n-1},...,d_{l_{n-1}}^{n-1})\right)
\]
with $d_{j}^i\in\N$ (we assume when using this notation that $d_1^i\geq d_2^i\geq...$).  Let $\mu^i=(d_{1}^i,...,d_{l_i}^i)$ denote the partition indexed by $i$ and let $\mu^{tw}$ correspond to the $n$-tuple of \textit{twisted} partitions $(\emptyset,\mu^1,...,\mu^{n-1})$.  At times it will be convenient to write $\mu$ as a multiset $\{\xi^id_{j}^i\}$ where the power of $\xi$ keeps track of which $\mu^i$ the $d_{j}^i$ came from.  Let $l(\mu):=\sum l_i$ denote the length of $\mu$.  Set $|\mu^i|:=\sum_j d_j^i$ and $|\mu|:=\sum |\mu^i|$.  Let $\underline\mu$ denote the underlying partition of $\mu$ that forgets the $\Z_n$ decorations.  We define $-\mu:=\{\xi^{-i}d_j^i\}$, i.e. it is the $n$-tuple of partitions with opposite twistings.  We also define
\[
z_\mu:=|\aut(\mu)|\prod nd_{j}^i
\]
to be the order of the centralizer of any element in the conjugacy class $\mu$.

Suppose $\lambda=(\lambda_0,...,\lambda_{n-1})$.  Via $n$-quotients (described explicitly in Section \ref{sec:wrfock}) $\lambda$ can be identified with a partition of $nd$ where $d=|\lambda|$.  We denote this corresponding partition by $\bar\lambda$.  We write $\bar\lambda=\{(i,j)\}$ where $i$ indexes the rows and $j$ indexes the columns of the Young diagram corresponding to $\bar\lambda$.  We will often think of $\bar\lambda$ as a colored Young diagram where the box $(i,j)$ has color $j-i \mod n$.  We denote the boxes with color $k$ by $\bar\lambda[k]$.  For $\square\in\bar\lambda$, we let $h_k(\square)$ denote the number of color $k$ boxes in the hook defined by $\square$ and we define \[n_k(\bar\lambda):=\sum_i(i-1)(\text{\# of color k boxes in the $i$th row}).\]

We let $\gamma$ denote a tuple of \textit{nontrivial} elements in $\Z_n$.  We define $m_i(\gamma)$ to be the number of occurrences of $\xi^i\in\Z_n$ in $\gamma$.

\subsection{Gromov-Witten Theory}

Given $\mu$ and $\gamma$ as above, let $\M_{g,\gamma+\mu}(\B\Z_n)$ denote the moduli stack of stable maps to the classifying space with $m_i(\gamma)+l_i(\mu)$ marked points twisted by $\xi^i$.  We recall the definitions of some natural classes on this moduli stack.

By the definition of $\B\Z_n$, $\M_{g,\gamma+\mu}(\B\Z_n)$ parametrizes degree $n$ covers of the source curve, ramified over the twisted points, with an action of $\Z_n$ which exhibits the source curve as a quotient of the cover.  Let
\[
p:\mathcal{U}_h\rightarrow \M_{g,\gamma+\mu}(\B\Z_n)
\]
be the universal covering curve of genus $h$ where $h$ is computed via the Riemann-Hurwitz formula.  The Hodge bundle on $\M_{g,\gamma+\mu}(\B\Z_n)$ is the rank $h$ bundle defined by
\[
\mathbb{E}:=p_*\omega_p
\]
where $\omega_p$ is the relative dualizing sheaf of $p$.  $\Z_n$ naturally acts on $\mathbb{E}$ and its dual $\mathbb{E}^\vee$.  For any $\zeta\in\Z_n$, we define $\mathbb{E}_\zeta$ and $\mathbb{E}_\zeta^\vee$ to be the $\zeta$-eigenbundles of $\mathbb{E}$ and $\mathbb{E}^\vee$, respectively.  They are related by the formula $(\mathbb{E}_\zeta)^\vee=\mathbb{E}_{\zeta^{-1}}^\vee$.  We also have the formula
\[
\mathbb{E}_{\zeta^{-1}}^\vee=R^1\pi_*f^*\so_{\zeta}
\]
where $\pi$ is the map from the universal curve, $f$ is the universal map, and $\so_{\zeta}$ is the line bundle with isotropy acting by multiplication by $\zeta$.  The lambda classes are defined as the chern classes of these bundles:
\[
\lambda_j^{\zeta}:=c_j\left(\mathbb{E}_\zeta\right)
\]

By forgetting the orbifold structure of the curve, there is a universal coarse curve
\[
q:\mathcal{U}_{g,|\gamma|+|\mu|}\rightarrow\M_{g,\gamma+\mu}(\B\Z_n)
\]
along with a section $s_p$ for each marked point $p$.  We define the cotangent line bundles by
\[
\mathbb{L}_p:=s_p^*\omega_q
\]
where $\omega_q$ is the relative dualizing sheaf of $q$. The psi classes on $\M_{g,\gamma+\mu}(\B\Z_n)$ are defined by
\[
\psi_p:=c_1\left(\mathbb{L}_p\right)
\]
The marked points in $\mu$ are indexed by $\{(i,j):0\leq i<n,1\leq j \leq l_i \}$ and we denote the corresponding psi classes by $\psi_{i,j}$.

For any $a\in\frac{1}{n}\Z$, the special cubic Hodge integrals we are interested in are
\begin{align}\label{eqn:vertex}
V_{g,\gamma}(\mu;a):=\nonumber\frac{(a+1)^{l_0}}{|\aut(\mu)|}&\prod_{i=0}^{n-1}\prod_{j=1}^{l_i}\frac{\prod_{k=0}^{d_{j}^i-1}(ad_{j}^i+\frac{i}{n}+k)}{(-1)^{d_j^i}d_{j}^i\cdot d_{j}^i!}\\
&\int_{\M_{g,\gamma+\mu}(\B\Z_n)}\frac{\Lambda^0(1)\Lambda^1(a)\Lambda^{-1}(-a-1)}{\delta(a)\prod_{i=0}^{n-1}\prod_{j=1}^{l_i}\left( \frac{1}{d_{j}^i}-\psi_{i,j} \right)}
\end{align}
where
\[
\Lambda^i(t):=(-1)^{\text{rk}}\sum_{j=0}^{\text{rk}}(-t)^{\text{rk}-j}\lambda_j^{\xi^i}
\]
with $\text{rk}:=\text{rk}\left(\mathbb{E}_{\xi^i}\right)$ and $\delta(a)$ the function which takes value $-a^2-a$ on the connected component of the moduli space which parametrizes trivial covers of the source and takes value $1$ on all other components.

\begin{remark}
The parameter $a$ is often referred to as the \textit{framing}.
\end{remark}


Introduce formal variables, $u$ and $x_i$ to track genus and marks.  Also introduce the variables $p_\mu$ with formal multiplication defined by concatenating the indexing partitions.  Then we define
\[
V_\mu^\bullet(x,u;a):=\exp\left(\sum_{g,\gamma,\nu} V_{g,\gamma}(\nu;a)u^{2g-2+l(\nu)}\frac{x^\gamma}{\gamma!}p_\nu\right)[p_\mu]
\]
where \[\frac{x^\gamma}{\gamma!}:=\prod_{i=1}^{n-1}\frac{x_i^{m_i(\gamma)}}{m_i(\gamma)!}\] and $[p_\mu]$ denotes ``the coefficient of $p_\mu$''.  By definition, $V_\mu^\bullet(x,u;a)$ is the one-leg $A_{n-1}$ orbifold GW vertex defined in \cite{r:lgoa}.

\begin{definition}
The framed GW vertex is defined by
\begin{equation}\label{modgw}
\tilde{V}_\mu^\bullet(a):=\prod_{i=1}^{n}(\sqrt{-1}\xi_{2n}^i)^{l_i}V_\mu^\bullet(x,u;a).
\end{equation}
where $l_n:=l_0$.
\end{definition}

\subsection{Donaldson-Thomas Theory}

Let $q_0,...,q_{n-1}$ be formal variables (always assume that the index of $q_k$ is computed modulo $n$) and define $q:=q_0\cdots q_{n-1}$.  For $\bar\lambda$ as above, define
\begin{equation}\label{dtvert}
P_{\lambda}(q_0,...,q_{n-1}):=\frac{1}{\prod_{\square\in\bar\lambda}\left(1-\prod_i q_i^{h_i(\square)}\right)}.
\end{equation}

By Theorem 12 in \cite{bcy:otv}, $P_{\lambda}(-q_0,...,q_{n-1})$ is the reduced, multi-regular one-leg $A_{n-1}$ orbifold DT vertex.

\begin{remark}\label{rmk:signdisc}
Notice the sign discrepancy between \eqref{dtvert} and the DT vertex.
\end{remark}

\begin{definition}
The framed DT vertex is defined by
\begin{equation}\label{moddt}
\hspace{-.75cm}\tilde{P}_\lambda(a):=\left(\left(\left(-\xi_{2n}\right)^{|\lambda|}\prod\xi_n^{l|\lambda_l|}\right)^n\prod_{(i,j)\in\bar\lambda}q_{j-i}^{j-i}\right)^{-a}\frac{\cha_{\bar\lambda}(n^d)}{\dim(\lambda)}q^{\frac{d}{2}}(-1)^d\prod_i q_i^{n_i(\bar\lambda)}P_{\lambda}(q_0,...,q_{n-1}).
\end{equation}
\end{definition}

\begin{remark}
$\cha_{\bar\lambda}$ is a character of $S_{dn}$ whereas $\dim(\lambda)$ is the dimension of an irreducible representation of $\Z_n\wr S_d$.  As we will see in Section \ref{sec:signs}, the quotient $\frac{\cha_{\bar\lambda}(n^d)}{\dim(\lambda)}$ is simply a compact way of keeping track of a sign.
\end{remark}

\begin{remark}
In Corollary \ref{sp}, we relate $\tilde{P}_\lambda(0)$ to loop Schur functions.
\end{remark}

\subsection{The Correspondence}

We will prove the following formula.

\begin{theorem1}
After the change of variables \[q\rightarrow \re^{\sqrt{-1}u},\hspace{.5cm} q_k\rightarrow \xi_n^{-1}\re^{-\sum_i\frac{\xi_n^{-ik}}{n}\left(\xi_{2n}^i-\xi_{2n}^{-i}\right)x_i} \hspace{.2cm}(k>0),\]
\[
\tilde{V}_\mu^\bullet(a)=\sum_\lambda \tilde{P}_\lambda(a)\frac{\cha_\lambda(\mu)}{z_\mu}
\]
\end{theorem1}

In Section \ref{sec:app}, we use Theorem \ref{thm1} to deduce the Gromov-Witten/Donaldson-Thomas correspondence for local $\Z_n$-gerbes over $\proj^1$.

\begin{theorem2}
Let $\X$ be a local $\Z_n$-gerbe over $\proj^1$ and let $GW(\X)$ and $DT_{mr}'(\X)$ denote the GW potential and the reduced, multi-regular DT potential of $\X$, respectively.  After the change of variables \[q\rightarrow -\re^{\sqrt{-1}u},\hspace{.5cm} q_k\rightarrow \xi_n^{-1}\re^{-\sum_i\frac{\xi_n^{-ik}}{n}\left(\xi_{2n}^i-\xi_{2n}^{-i}\right)x_i} \hspace{.2cm}(k>0),\]
\[
GW(\X)=DT_{mr}'(\X).
\]
\end{theorem2}

\begin{remark}
Notice the sign difference in the change of variables of Theorems 1 and 2 -- this difference is an artifact of Remark \ref{rmk:signdisc}.
\end{remark}

\begin{remark}
The change of variables in Theorems 1 and 2 is predicted by Iritani's stacky Mukai vector \cite{i:rcaisiqc} and previously appeared in \cite{z:gmvf}.  We thank Jim Bryan for explaining this change of variables to us.
\end{remark}

\section{Geometry}\label{sec:geometry}

In this section we set up auxilary integrals on moduli spaces of relative maps into $\proj^1$-gerbes in order to obtain bilinear relations between the vertex $\tilde V_\mu^\bullet(a)$ and certain rubber integrals $\tilde H_{\nu,\mu}^\bullet(a)$.  The rubber integrals in $\tilde H_{\nu,\mu}^\bullet(a)$ can be interpreted as wreath Hurwitz numbers and can be computed via Burnside's formula in terms of the represenation theory of the wreath product $\Z_n\wr S_d$.  We use this interpretation in Sections \ref{sec:invertible} and \ref{sec:linearalgebra}  to show that the localization relations uniquely determine $\tilde V_\mu^\bullet(a)$ from  $\tilde H_{\nu,\mu}^\bullet(a)$.  The method of localizing maps into gerbes in order to obtain useful relations of Hodge integrals first appeared in work of Cadman and Cavalieri \cite{cc:gl} where it was used to compute the GW invariants of $[\C^3/\Z_3]$.

\subsection{Cyclic Gerbes over $\proj^1$}\label{sec:gerbes}

Cyclic $\proj^1$ gerbes will be important both for the localization computations in Section \ref{sec:bil} and in the GW/DT comparisons in Section \ref{sec:app}.  We briefly collect the necessary details here.  For each line bundle $\so(-k)$ with $0\leq k<n$, we can define a $\proj^1$-gerbe $\G_k$ with isotropy group $\Z_n$ and an orbifold line bundle $L_k$ as follows.

\begin{definition}\label{gerbe}
The gerbe $\G_k$ is defined by pullback
\[\begin{CD}
\G_k   	@>>>          \B\C^*       \\
@VVV                             @VV\lambda\rightarrow\lambda^nV\\
\proj^1                  @>\so(-k)>>      \B\C^*
\end{CD}\]
and $L_k$ is defined to be the line bundle parametrized by the top map.
\end{definition}

Note that the numerical degree of $L_k$ is $-k/n$ and the action of $\Z_n$ on the fibers is given by multiplication by $\xi_n$ (see e.g. Section 2.3 of \cite{r:lgoa}).

The $\G_k$ are only distinct if we choose an isomorphism of each isotropy group with $\Z_n$.  In other words, for each $\phi\in\aut(\Z_n)$, we obtain an equivalence $\tilde\phi_k:\G_k\stackrel{\cong}{\longrightarrow}\G_{\phi(k)}$ for each $k$.  However, it is not true in general that $\tilde\phi_k^*\left(L_{\phi(k)}\right)=L_k$.  This fact will be important in our discussion of 3-fold targets in Section \ref{sec:app}.

One of the most useful aspects of localizing maps of curves into $\proj^1$ gerbes is that it allows us to control the orbifold structure over $0$ and $\infty$.  To make this precise, let $\mathcal{C}$ be an orbifold with coarse space $\proj^1$ and orbifold structure only at $0$ and $\infty$.  Let $f:\mathcal{C}\rightarrow \G_k$ be a $\C^*$ fixed degree $d$ map with twisting $k_0$ at $0$ and $k_\infty$ at $\infty$.  Then
\[
k_\infty = -dk-k_0 \mod n.
\]
A more general characterization of this property was given in Lemmas II.12 and II.13 of Johnson's thesis \cite{j:egwtods}.  To keep track of this twisting compatibility, we make the following definition.

\begin{definition}\label{def:gk}
For a decorated partition $\mu=\{\xi^id_j^i\}$, we define the involution $g_k(\mu)$ by
\[
g_k(\mu):=\{\xi^{d_j^ik-i}d_j^i\}
\]
\end{definition}

If $f$ is a $\C^*$ fixed map from a disjoint union of orbifold $\proj^1$s with degree and twisting over $0$ given by $\mu$, then the degree and twisting over $\infty$ is determined by $-g_k(-\mu)$ (the conventions with signs seems cumbersome at the moment but it will be natural in later formulas).

\subsection{Auxilary Integrals}\label{sec:aux}

Here we set up integrals on the moduli spaces $\M_{g,\gamma}(\G_k,\mu[\infty])$ which parametrize maps with fixed ramification and isotropy profile over $\infty$.  These moduli spaces were developed by Abramovich and Fantechi \cite{af:otidf}.  The integrals we will investigate are the following.
\begin{equation}\tag{I-1}\label{int1}
\frac{1}{|\aut(\mu)|}\int_{\M_{g,\gamma}(\G_0,\mu[\infty])}\re(R^1\pi_*((\hat f^*L_0)(-D)\oplus \hat f^*L_0^\vee(-1)))
\end{equation}
where $D$ is the locus of relative points on the universal curve with trivial isotropy and $\hat f$ contracts the degenerated target and maps all the way to $\G_0$, and for $1\leq k\leq n-1$
\begin{equation}\tag{I-2}\label{int2}
\frac{1}{|\aut(\mu)|}\int_{\M_{g,\gamma}(\G_k,\mu[\infty])}\re(R^1\pi_*(\hat f^*L_k\oplus \hat f^*L_k^\vee(-1))).
\end{equation}

\subsection{Partial Evaluations}\label{sec:partialevals}

In certain cases, we can evaluate the integrals \eqref{int1} and \eqref{int2} explicitly.  We collect these computations here.

We begin with the first integral.  As we will see in Section \ref{sec:bil}, \eqref{int1} is equal to $V_{g,\gamma}(\mu;0)$.  Therefore, we consider special choices of $\mu$ for which we can evaluate $V_{g,\gamma}(\mu;0)$.  Recall that $\{d\}$ denotes the $n$-tuple of partitions with one untwisted part.  The following evaluation will be extremely useful.

\begin{lemma}\label{lemmag1}
\[
V_{g,\gamma}(\{d\};0)=\delta_{|\gamma|,0}\frac{(-1)^{d-1}}{n}\int_{\M_{g,1}}\lambda_g(d\psi)^{2g-2}.
\]
\end{lemma}

\begin{proof}
By \eqref{eqn:vertex}, $V_{g,\gamma}(\{d\};0)$ vanishes away from the locus of maps which parametrize trivial covers.  In particular, since $\gamma$ consists of nontrivial elements in $\Z_n$, the cover can only be trivial if $\gamma=\emptyset$.  On the locus of maps which parametrize trivial covers, $\mathbb{E}_\xi\cong\mathbb{E}_{\xi^{-1}}\cong\mathbb{E}_1$.  Therefore we can apply the Mumford relation to the integrand in the definition of $V_{g,\emptyset}((d);0)$.  The lemma follows by pushing forward to $\M_{g,1}$ which is a degree $\frac{1}{n}$ map.
\end{proof}

\begin{corollary}\label{corg2}
\[
V_{\mu}^\bullet(0)=\left(\frac{1}{z_{\mu^0}}\prod_{j=1}^{l_0}\frac{(-1)^{d_j^0-1}}{2}\csc\left( \frac{d_{j}^0u}{2}\right)\right) V_{\mu^{tw}}^\bullet(0)
\]
\end{corollary}

\begin{proof}
By \eqref{eqn:vertex}, the only nonzero vertex terms $V_{g,\gamma}(\mu,0)$ with $\mu^0\neq\emptyset$ are those with $l_0=1$ -- these invariants were computed in Lemma \ref{lemmag1}.  The evaluations of Lemma \ref{lemmag1} can be packaged using the Faber-Pandharipande identity \cite{fp:hiagwt}:
\[
\sum_g\left(\int_{\M_g,1}\lambda_g\psi^{2g-2}\right)t^{2g}=\frac{t}{2}\csc\left(\frac{t}{2}\right).
\]
The result then follows by passing from the connected invariants to the disconnected ones by exponentiating.
\end{proof}

From these evaluations, we see that the $a=0$ vertex is completely determined once we know the contributions coming from partitions $\mu$ with $\mu^0=\emptyset$.

For the integral \eqref{int2}, we obtain the following vanishing result.

\begin{lemma}
The integral \eqref{int2} vanishes if any of the parts of $\mu$ are untwisted.
\end{lemma}

\begin{proof}
The integral vanishes by dimensional reasons.  The dimension of the moduli space is $|\mu|+2g-2+|\gamma|+l(\mu)$.  The degree of the integrand is $|\mu|+2g-2+|\gamma|+l(\mu^{tw})$ which can be computed by the orbifold Riemann-Roch formula \cite[Theorem 7.2.1]{agv:gwtdms}.
\end{proof}

\subsection{Bilinear Relations}\label{sec:bil}

We now compute the integrals \eqref{int1} and \eqref{int2} via localization.  Beginning with \eqref{int1}, we give the target the standard $\C^*$ action with weight $1$ ($-1$) on the fibers of the tangent bundle over $0$ ($\infty$).  This defines a $\C^*$ action on the moduli space by postcomposing the map with the action.  In order to choose an equivariant lift of the integrand, we lift the action from the target to the bundles $T(-\infty)$, $L_0^\vee$, and $L_0(-1)$ so that $\C^*$ acts on the fibers over $0$ and $\infty$ with weights summarized in the following table.

\begin{center}
\begin{tabular}{c|ccc}
&$T(-\infty)$&$L_0$&$ L_0^\vee(-1)$\\
\hline
$0$&$1$&$a$&$-a-1$\\
$\infty$&$0$&$a$&$-a$
\end{tabular}
\end{center}

Each fixed locus of the torus action on the moduli space can be encoded by a bipartite graph $\Gamma$ with white (black) vertices corresponding to the connected components of $\hat f^{-1}(0)$ ($\hat f^{-1}(\infty)$).  The vertices and edges are decorated with the following data:
\begin{itemize}
\item Each vertex $v$ is labeled with a tuple $\gamma_v$ of nontrivial elements in $\Z_n$ corresponding to the twisted marks on that component and an integer $g_v$ corresponding to the genus.
\item Each edge $e$ is labeled with a complex number $(\xi^{k_e}d_e)$ which induces a $n$-tuple of partitions $\nu_v\in \text{Conj}(\Z_n\wr S_{d_v})$ at each white vertex and $-\nu_v\in \text{Conj}(\Z_n\wr S_{d_v})$ at each black vertex.
\item In addition, each black vertex is labeled with a $n$-tuple of partitions $\mu_v$ such that $|\mu_v|=|\nu_v|$ and the union of all $\mu_v$ is $\mu$.
\end{itemize}

To a white vertex, we associate the contribution
\[
\text{Cont}(v)=V_{g_v,\gamma_v}(\nu_v;a)
\]
and to a black vertex we associate the contribution
\begin{align*}
\text{Cont}(v)=&\frac{(-1)^{l_0(\nu_v)+g-1+\sum_{i\neq 0}\frac{n-i}{n}(m_i(\gamma_v)+l_i(\mu_v)+l_{n-i}(\nu_v))}(a)^{2g_v-2+|\gamma_v|+l(\mu_v)+l(\nu_v)}}{|\aut(\nu_v)|}\\
&\hspace{1cm}\cdot\left( \prod_{i=1}^{l(\nu_v)}nd_i \right)\int_{\M_{g_v,\gamma_v}(\G_0;-\nu_v[0],\mu_v[\infty])//\C^*}-(-\psi_0)^{2g_v-3+|\gamma_v|+l(\nu_v)+l(\mu_v)},
\end{align*}
where $\psi_0$ is the target psi class.  By the localization formula for orbifold stable maps (see for example \cite{cc:gl,l:ligwtaogwt,r:lgoa,z:gmvf}) we compute the integral
\[
\eqref{int1}=\frac{1}{|\aut(\mu)|}\sum_\Gamma \frac{1}{|\aut(\Gamma)|}\prod_v\text{Cont}(v).
\]

\begin{remark}
In the simplification of the black vertex contribution, we used the $\Z_n$-Mumford relation proved by Bryan, Graber and Pandharipande \cite{bgp:crc}, namely:
\[
\Lambda^1(a)\Lambda^{-1}(-a)=(a)^{rk(\mathbb{E}_{\xi})}(-a)^{rk(\mathbb{E}_{\xi^{-1}})}
\]
where the ranks can be computed by the orbifold Riemann-Roch formula.
\end{remark}


Setting $a=0$, we observe that the contributions from black vertices vanish and the integral is equal to $V_{g,\gamma}(\mu;0)$, justifying the statement made at the beginning of Section \ref{sec:partialevals}.

Define the rubber integral generating function
\[
H_{\nu,\mu}(x,u):=\frac{1}{|\aut(\nu)||\aut(\mu)|}\sum_{g,\gamma}\int_{\M}\psi_0^{r+|\gamma|-1}u^r\frac{x^\gamma}{\gamma!}
\]
where $r:=2g-2+l(\mu)+l(\nu)$, $\M$ is the space of relative maps into the rubber target: ${\M_{g,\gamma}(\G_0;\nu[0],\mu[\infty])//\C^*}$.

For notational convenience, we define
\[
\tilde{H}_{\nu,\mu}^\bullet(a):=\exp\left( H_{\nu,\mu}(a\xi_{2n}^{-1}x_1,...,a\xi_{2n}^{1-n}x_{n-1},\sqrt{-1}au) \right)
\]

The above localization computations amount to the following bilinear relations between $V$ and $H$:
\begin{equation}\tag{R-1}\label{rel1}
\tilde{V}_\mu^\bullet(0)=\sum_{|\nu|=|\mu|}\tilde{V}_\nu^\bullet(a)z_\nu\tilde{H}_{-\nu,\mu}^\bullet(a).
\end{equation}

\begin{remark}
Notice that the $-\nu$ appearing in the rubber integrals is equal to $g_0(\nu)$ from Definition \ref{def:gk}.
\end{remark}


We also compute \eqref{int2} via localization.  Again we equip the moduli space with a $\C^*$ action via the standard $\C^*$ action on the target.  We lift the integrand with the choice of linearizations summarized in the following table.

\begin{center}
\begin{tabular}{c|ccc}
&$T(-\infty)$&$L_k$&$ L_k^\vee(-1)$\\
\hline
$0$&$1$&$0$&$-1$\\
$\infty$&$0$&$k/n$&$-k/n$
\end{tabular}
\end{center}

The localization computation of \eqref{int2} is nearly identical to that of \eqref{int1} and leads to the relations
\begin{equation}\tag{R-2}\label{rel2}
0=\sum_{|\nu|=|\mu|}\tilde{V}_\nu^\bullet(0)z_\nu\tilde{H}_{g_k(\nu),\mu}^\bullet\left(\frac{k}{n}\right)
\end{equation}
where $\mu$ is any partition with at least one untwisted part.


\subsection{Wreath Hurwitz Numbers}\label{sec:wreath}

In the non-orbifold case, it was shown by Liu, Liu, and Zhou \cite{llz:pmvf,llz:ftphi} that certain rubber integrals can be interpreted in terms of double Hurwitz numbers.  In this section, we generalize their result to the orbifold case.

Hurwitz numbers classically count degree $d$ ramified covers of Riemann surfaces with monodromy around the branch points prescribed by conjugacy classes in $S_d$.  Cyclic wreath Hurwitz numbers are defined to be analogous counts of degree $dn$ ramified covers where the monodromy is prescribed by conjugacy classes $\mu$ in $\Z_n\wr S_d$.  Since $\Z_n$ is in the center of $\Z_n\wr S_d$, such covers have a natural $\Z_n$ action and the quotient is a classical Hurwitz cover with monodromy given by the underlying partitions $\underline\mu$.

We define now the particular wreath Hurwitz numbers which arise in our context.

\begin{definition}
Let $H_{\nu,\mu}^{g,\gamma}$ be the automorphism-weighted count of wreath Hurwitz covers $f:C\rightarrow\proj^1$ where the branch locus consists of a set of $|\gamma|+r+2$ fixed points (we fix the last two points at $0$ and $\infty$) and the maps satisfy the following conditions:
\begin{itemize}
\item The quotient $C/\Z_n$ is a connected genus $g$ curve,
\item The monodromy around $0$ and $\infty$ is given by $\nu$ and $\mu$,
\item The monodromy around the branch point corresponding to $\gamma_i\in\gamma$ is given by the conjugacy class $\{\gamma_i,1,...,1\}$,
\item The monodromy around the $r$ additional branch points is given by the conjugacy class $\{2,1,...,1\}$.
\end{itemize}
\end{definition}

\begin{remark}
Here we use the multiset notation for $n$-tuples of partitions introduced in Section \ref{sec:parts}.
\end{remark}

The next theorem relates the rubber integrals which arose in the localization computations to the wreath Hurwitz numbers $H_{\nu,\mu}^{g,\gamma}$.

\begin{theorem}\label{thm:whn}
\[
H_{\nu,\mu}^{g,\gamma}=\frac{r!}{|\aut(\nu)||\aut(\mu)|}\int_{\M_{g,\gamma}(\G_0;\nu[0],\mu[\infty])//\C^*}\psi_0^{r-1+|\gamma|}.
\]
\end{theorem}

\begin{proof}
Via the forgetful map $F:\M_{g,\gamma}(\G_0;\nu[0],\mu[\infty])\rightarrow\M_{g,n}(\proj^1;\underline\nu[0],\underline\mu[\infty])$, we obtain a branch morphism $Br:\M_{g,\gamma}(\G_0;\nu[0],\mu[\infty])\rightarrow \text{Sym}^r\proj^1\cong\proj^r$ by postcomposing $F$ with the usual branch morphism.  For each of the $n$ (twisted) marked points, we also obtain maps $\tilde{ev}_i:\M_{g,\gamma}(\G_0;\nu[0],\mu[\infty])\rightarrow \proj^1$ by postcomposing the usual evaluation map with the natural map to $\proj^1$.  Then the wreath Hurwitz numbers can be expressed as
\begin{equation}\label{hurwitzint}
H_{\nu,\mu}^{g,\gamma}=\frac{1}{|\aut(\nu)||\aut(\mu)|}\int_{\M_{g,\gamma}(\G_0;\nu[0],\mu[\infty])}Br^*(pt)\cdot\prod \tilde{ev}_i^*(pt).
\end{equation}

It is left to show that
\[
\int_{\M_{g,\gamma}(\G_0;\nu[0],\mu[\infty])}Br^*(pt)\cdot\prod \tilde{ev}_i^*(pt)=r!\int_{\M_{g,\gamma}(\G_0;\nu[0],\mu[\infty])//\C^*}\psi_0^{r-1+|\gamma|}
\]
and we accomplish this via localization.

We equip the moduli space with a torus action by fixing the $\C^*$ action on the target $t\cdot[z_0:z_1]=[z_0:tz_1]$ so that the tangent bundle is linearized with weights $1$ at $0=[0:1]$ and $-1$ and $\infty=[1:0]$.  The isomorphism $\proj^r=\proj(H^0(\proj^1,\so(r)))\rightarrow \text{Sym}^r\proj^1$ is given by $s\rightarrow Div(s)$ where the basis $\langle z_0^r,z_0^{r-1}z_1,...,z_1^r \rangle$ for $H^0(\proj^1,\so(r))$ corresponds to the homogeneous coordinates $(y_0:y_1:...:y_r)$.  We equip $\proj^r$ with the torus action $t\cdot(y_0:y_1:...:y_r)=(y_0:ty_1:...:t^ry_r)$ which makes $Br$ an equivariant map.  We lift $[pt]\in H^{2r}(\proj^r)$ to $\prod_{i=0}^{r-1}(H+i\hbar)\in H_{\C^*}^{2r}(\proj^r)$ where $\hbar$ is the equivariant parameter.  The preimage of this lift is the locus of maps where the simple ramification points map to $\infty$.  Likewise we lift
\[
\tilde{ev}_i^*(pt)=c_1(\tilde{ev}_i^*\so(1))
\]
by linearizing $\so(1)$ with weights $0$ at $0$ and $-1$ at $\infty$.

With these choices of linearizations, we see that the integrand vanishes on all fixed loci where any of the $n+r$ points with nontrivial monodromy map to $0$.  This leaves exactly one fixed locus where the target expands over $\infty$ and everything interesting happens over the expansion.  On this locus, the integrand specializes to $(-\hbar)^{r+n}r!$ and the inverse of the equivariant Euler class of the normal bundle is
\[
\frac{1}{-\hbar-\psi_0}.
\]
Therefore the contribution, and hence the integral in \eqref{hurwitzint}, is equal to
\[
r!\int_{\M_{g,\gamma}(\G_0;\nu[0],\mu[\infty])//\C^*}\psi_0^{r+n-1}.
\]
\end{proof}

\begin{corollary}
\begin{align*}
H^\bullet_{\nu,\mu}(x,u)&=\exp\left( \sum_{g,\gamma} H_{\nu,\mu}^{g,\gamma}\frac{u^r}{r!}\frac{x^\gamma}{\gamma!} \right)\\
&=\sum_{g,\gamma} H_{\nu,\mu}^{\cha,\gamma \bullet}\frac{u^r}{r!}\frac{x^\gamma}{\gamma!}
\end{align*}
where $H_{\nu,\mu}^{\cha,\gamma \bullet}$ is the wreath Hurwitz number with possibly disconnected covers.
\end{corollary}

By the Burnside formula (see Dijkgraaf \cite{d:msaec} for a proof), we compute
\[
H_{\nu,\mu}^{\cha,\gamma \bullet}=\sum_{|\lambda|=d}\left(f_T(\lambda)\right)^r \prod\left(f_{i}(\lambda)\right)^{m_i(\gamma)}\frac{\cha_{\lambda}(\mu)}{z_{\mu}}\frac{\cha_{\lambda}(\nu)}{z_{\nu}}
\]
where $f_T(\lambda)$ and $f_i(\lambda)$ are the \textit{central characters} defined by
\[
f_T(\lambda):=\frac{nd(d-1)\cha_{\lambda}(\{2,1,...,1\})}{2\cdot\text{dim}\lambda}
\]
and
\[
f_{i}(\lambda):=\frac{d\cha_{\lambda}(\{\xi^i,1,...,1\})}{\text{dim}\lambda}.
\]

Therefore we obtain the following form for the generating function of wreath Hurwitz numbers:
\begin{equation}\label{hurgen}
H_{\nu,\mu}^\bullet(x,u)=\sum_{|\lambda|=d}\frac{\cha_{\lambda}(\mu)}{z_{\mu}}\frac{\cha_{\lambda}(\nu)}{z_{\nu}}\re^{f_T(\lambda)u+\sum f_{i}(\lambda)x_i}.
\end{equation}

Using the fact that $\cha_\lambda(-\nu)=\overline{\cha_\lambda(\nu)}$, orthogonality of characters gives us the following relations:
\begin{equation}\label{hurwitz1}
H_{\nu,\mu}^\bullet(x+y,u+v)=\sum_\sigma H_{\nu,\sigma}^\bullet(x,u)z_\sigma H_{-\sigma,\mu}^\bullet(y,v)
\end{equation}
and
\begin{equation}\label{hurwitz2}
H_{\nu,-\mu}^\bullet(0,0)=\frac{1}{z_\mu}\delta_{\nu,\mu}.
\end{equation}

The relations \eqref{hurwitz1} and \eqref{hurwitz2} also have a geometric meaning -- \eqref{hurwitz1} is the degeneration formula for the target $\proj^1$ where $x$ and $y$ keep track of whether the corresponding point of ramification maps to one side of the node or the other, and \eqref{hurwitz2} counts covers with ramification only over $0$ and $\infty$.

\subsection{Invertibility}\label{sec:invertible}

In this section we show that the relations \eqref{rel1} can be inverted explicitly.  We also state the main result concerning the relations \eqref{rel2} but we defer the proof to the next section.

The next lemma follows immediately from Equations \eqref{hurwitz1} and \eqref{hurwitz2}.

\begin{lemma}\label{framing1}
Framing dependence in the conjugacy basis:
\[
\tilde{V}_\mu^\bullet(a)=\sum_{|\nu|=|\mu|}\tilde{V}_\nu^\bullet(0)z_\nu\tilde{H}_{-\nu,\mu}^\bullet(-a)
\]
\end{lemma}

In particular, Lemma \ref{framing1} determines the general framed vertex from the $a=0$ vertex and characters of $\Z_n\wr S_d$.

Define
\[
\hat{P}_\lambda(a):=\sum_\mu \tilde{V}_\mu^\bullet(a)\cha_\lambda(-\mu)
\]
or equivalently
\[
\tilde{V}_\mu^\bullet(a)=\sum_\lambda \hat{P}_\lambda(a)\frac{\cha_\lambda(\mu)}{z_\mu}.
\]

Then Lemma \ref{framing1} is equivalent to the following.

\begin{lemma}\label{framing2}
Framing dependence in the representation basis:
\[
\hat{P}_\lambda(a)=\re^{-\sqrt{-1}af_T(\lambda)u-a\sum \xi_{2n}^{-i} f_{i}(\lambda)x_i}\hat{P}_\lambda(0)
\]
\end{lemma}

Therefore, once we know that $\tilde{P}_\lambda(a)$ and $\tilde{P}_\lambda(0)$ are related by the exponential factor of Lemma \ref{framing2}, we only need to prove Theorem \ref{thm1} for the case $a=0$.

The relations \eqref{rel2} are significantly more difficult to work with and do not admit a convenient inverse as far as we know.  Nonetheless, we prove that they are invertible.

\begin{theorem}\label{thm:invert}
Relations \eqref{rel2} uniquely determine $V_\mu(0)$ from the partial evaluations of Corollary \ref{corg2} and characters of $\Z_n\wr S_d$.
\end{theorem}

The proof of Theorem \ref{thm:invert} is rather involved and we defer it to the next section.  In the meantime, we gather formally the reductions which we have made while the formulas are fresh in our minds.

\begin{reduction}\label{reduction}
To prove Theorem \ref{thm1}, it suffices to check that the following properties hold after the prescribed change of variables.
\begin{enumerate}[(I)]
\item The framing factors are consistent:
\[
\left(\left(\left(-\xi_{2n}\right)^{|\lambda|}\prod\xi_n^{l|\lambda_l|}\right)^n\prod_{(i,j)\in\bar\lambda}q_{j-i}^{j-i}\right)^a=\re^{\sqrt{-1}af_T(\lambda)u+a\sum \xi_{2n}^{-i} f_{i}(\lambda)x_i}.
\]

\item $\tilde{P}_\lambda(0)$ satisfy the partial evaluations of Corollary \ref{corg2}:
\[
\sum_{|\lambda|=|\mu|} \tilde{P}_\lambda(0)\frac{\cha_\lambda(\mu)}{z_\mu}=\left(\frac{1}{z_{\mu^0}}\prod_{j=1}^{l_0}\frac{\sqrt{-1}(-1)^{d_j^0}}{2}\csc\left( \frac{d_{j}^0u}{2}\right)\right)\left( \sum_{|\sigma|=|\mu^{tw}|} \tilde{P}_{\sigma}(0)\frac{\cha_{\sigma}(\mu^{tw})}{z_{\mu^{tw}}} \right).
\]

\item $\tilde{P}_\lambda(0)$ satisfy the relations \eqref{rel2} for all $\mu$ with at least one untwisted part:
\[
\sum_\nu\left( \sum_\lambda \tilde{P}_\lambda(0)\frac{\cha_\lambda(\nu)}{z_\nu} \right)z_\nu \left( \sum_\sigma \frac{\cha_\sigma(g_k(\nu))}{z_{g_k(\nu)}}\frac{\cha_\sigma(\mu)}{z_\mu}\re^{\frac{k}{n}(\sqrt{-1}f_T(\sigma)u+\sum \xi_{2n}^{-i}f_{i}(\sigma)x_i)}\right)=0.
\]

\end{enumerate}
\end{reduction}

\begin{proof}
If $\tilde P_\lambda(0)$ satisfies (II) and (III), then Corollary \ref{corg2} and Theorem \ref{thm:invert} imply that Theorem \ref{thm1} is true in the case $a=0$.  The general framed correspondence then follows from the definition of the framed DT vertex and Lemma \ref{framing2}.
\end{proof}

The proofs of identities (I) - (III) are given in Section \ref{sec:proof} after developing the necessary combinatorial and representation theoretic identities in Sections \ref{sec:combinatorics} and \ref{sec:reptheory}.

\section{Linear Algebra}\label{sec:linearalgebra}

This section is devoted to the proof of Theorem \ref{thm:invert}.  By Corollary \ref{corg2}, the only vertices left to be determined are those $V^\bullet_{\nu}(0)$ with $\nu^{tw}=\nu$. So let us rewrite \eqref{rel2} as
\begin{equation}\tag{R-2'}\label{rel2'}
0=\sum_{|\nu|=|\mu|}\tilde{V}_{\nu^{tw}}^\bullet(0)\tilde{V}_{\nu^0}^\bullet(0)
z_\nu\tilde{H}_{g_k(\nu),\mu}^\bullet\left(\frac{k}{n}\right).
\end{equation}

Let us begin by reinterpreting relations \eqref{rel2'} in terms of matrix equations.  Define the vector
\[
\alpha_d=(\tilde{V}_{\nu}^\bullet(0))
\]
with indexing set $\{\nu:|\nu|\leq d,\nu=\nu^{tw}\}$ and the vector
\[
\beta_d=\left(-\sum_{|\tau|=|\mu|,\tau=\tau^0}\tilde{V}_{\tau}^\bullet(0)
z_\tau\tilde{H}_{g_k(\tau),\mu}^\bullet\left(\frac{k}{n}\right)\right)
\]
with indexing set $\{(\mu,k):|\mu|\leq d,\mu\neq\mu^{tw}, k\neq 0\}$.  We introduce a matrix
\[
\Phi_d(u;x)=(\Phi_d^{(\mu,k),\nu}(u;x))_{(\mu,k),\nu}
\]
with the same indexing sets defined by
\[
\Phi_d^{(\mu,k),\nu}(u;x)=\left\{\begin{array}{ll}0,&\textrm{if}|\nu|>|\mu|\\
z_\nu\tilde{H}_{g_k(\nu),\mu}^\bullet\left(\frac{k}{n}\right),&\textrm{if}|\nu|=|\mu|\\
\sum_{|\tau|=|\mu|-|\nu|,\tau^0=\tau}\tilde{V}_{\tau}^\bullet(0)z_{(\tau\sqcup\nu)}
\tilde{H}_{g_k(\tau\sqcup\nu),\mu}^\bullet\left(\frac{k}{n}\right),
&\textrm{if}|\nu|<|\mu|\end{array} \right.
\]

Then the collection of relations \eqref{rel2'} is equivalent to the collection of matrix equations
\[
\Phi_d(u;x)\alpha_d=\beta_d
\]
Our task is to show that $\Phi_d(u;x)$ has full (column) rank for all $d$.  

\subsection{Matrix Reductions}
We begin by making a sequence of reductions.  First note that
\[\Phi_d=\left(
\begin{array}{cccc}
\Phi'_d&*&\cdots&*\\
0&\Phi'_{d-1}&\ddots&\vdots\\
\vdots&\ddots&\ddots&*\\
0&\cdots&0&\Phi'_1
\end{array} \right)
\]
where $\Phi_d'$ is defined by restricting the indexing sets to partitions of size $d$.  Therefore, it suffices to prove that $\Phi_d'$ has full rank and to do this we need only prove that the specialization $\tilde{\Phi}_d:=\Phi_d'|_{x_2=...=x_{n-1}=u=0}$ has full rank.

\begin{remark}\label{clarify}
By setting $u=0$, notice that $\tilde{\Phi}_d$ is a generating function of wreath Hurwitz numbers counting covers for which the $\Z_n$ quotient is a disjoint union of $\proj^1$s, each one fully ramified over $0$ and $\infty$.  Moreover, by setting $x_2=\dots = x_{n-1}=0$, the only nontrivial monodromy away from $0$ and $\infty$ is given by conjugacy classes $\{\xi,1,\dots,1\}$.
\end{remark}

By the first part of Remark \ref{clarify}, if $\underline\mu\neq \underline\nu$, then the entry $\tilde\Phi_d^{(\mu,k),\nu}=0$.  This implies
\[
\tilde\Phi_d=\left(
\begin{array}{cccc}
\tilde\Phi_{\tau^1}&0&\cdots&0\\
0&\tilde\Phi_{\tau^2}&\ddots&\vdots\\
\vdots&\ddots&\ddots&0\\
0&\cdots&0&\ddots
\end{array} \right)
\]
where $\tilde\Phi_\tau$ is defined by restricting the indexing sets to a single underlying partition $\underline\nu=\underline\mu=\tau$ of size $d$.  Therefore, we have reduced our task to showing that $\tilde\Phi_\tau$ has full rank for a fixed partition $\tau$.

To this end, fix $\tau$ once and for all and write $\tau=(\tau_1,...,\tau_l)$ with nonincreasing parts.  We henceforth suppress $\tau$ from the notation and write $\tilde\Phi$ for $\tilde\Phi_\tau$.  We also write $x$ for $x_1$ when no confusion arises.

In order to prove that $\tilde\Phi$ has full rank, we will restrict the row index to a suitable subset and show that the resulting submatrix is invertible.  In order to do this, we must first introduce some subtle notation.

\subsection{Ordering Convention}

We introduce an order on the parts of each $\mu$ with $\underline{\mu}=\tau$.  Begin by defining $c_i:=\Gcd(\tau_i,n)$. If $\underline\mu=\tau$, we can write $\mu$ as the multiset $\mu=\{\xi^{t_i}\tau_i\}$ where $t_i\in\{0,...,n-1\}$.  We define $\bar t_i:=t_i($mod $c_i)$ and we set $(\tau_i,t_i)>(\tau_j,t_j)$ if one of the following is true
\begin{enumerate}
\item $\tau_i>\tau_j$, or
\item $\tau_i=\tau_j$ and $\bar t_i<\bar t_j$, or
\item $\tau_i=\tau_j$ and $\bar t_i=\bar t_j$ and $t_i<t_j$.
\end{enumerate}
Then $\mu$ can be written uniquely as
\[
\mu=((\tau_1,m_1),\dots,(\tau_l,m_l))
\]
where the pairs are nonincreasing.  This ordering convention will be important in defining the square submatrix $\hat\Phi$ in Section \ref{subsec:square} and in proving its invertibility in Section \ref{subsec:invert}.  At present, we use the ordering convention to define \[\tilde{\mu}:=\mu\setminus(\tau_1,m_1)\] and we define the \textit{twisting partition} of $\mu$ to be \[t(\mu):=(m_1,...,m_l).\]

\subsection{A Square Submatrix}\label{subsec:square}

We now explain a particular way to reduce the row index of $\tilde\Phi$ to a suitable subset so that the resulting submatrix $\hat\Phi$ is square.  We will show in Section \ref{subsec:invert} that $\hat\Phi$ is invertible which proves that $\tilde{\Phi}$ has full rank.

For $d\geq 1$ and $h\in\{1,...,n-1\}$, let $c:=\Gcd(n, d)$ and $\bar h=h($mod $c)\in\{0,...,n-1\}$.  We define
\[
\Sigma_{d,h}:=\{k\in\{1,\cdots,n-1\}|-h+dk=-\bar{h} (\text{mod } n)\}
\]

\begin{remark}
$\bar h$ has the following interpretation: For each $k\in\{1,...,n-1\}$ consider the unique $\C^*$ fixed map from an effective orbifold $\proj^1$ with orbifold ramification $(d,h)$ at $0$. Then the twisting over $\infty$ is fixed (c.f. Section \ref{sec:gerbes}) and $\bar h$ is the smallest possible twisting at $\infty$ as we vary $k$.  Moreover, $\Sigma_{d,h}$ is exactly the set of $k$ for which the minimal twisting is obtained.
\end{remark}

Notice that
\[
|\Sigma_{d,h}|=\begin{cases}
c-1 &\text{ if } h\in\{1,\dots,c-1\}\\
c &\text{ else. }
\end{cases}
\]
The set $\Sigma_{d,h}$ has a natural order as a subset of $\{1,...,n-1\}$, so we can write $\Sigma_{d,h}=\{k_1,...,k_{|\Sigma_{d,h}|}\}$.  We define
\[
k_d(h):=\begin{cases}
k_{\bar h}\in\Sigma_{d,h} &\text{ if } h\in\{1,...,c-1\}\\
k_{\bar h+1}\in\Sigma_{d,h} &\text{ else.}
\end{cases}
\]

\begin{lemma}\label{lem:bij}
$k_d(-)$ defines a bijection on the set $\{1,...,n-1\}$.
\end{lemma}
\begin{proof}
We show that the map is injective.  Suppose $k_d(h)=k_d(h')$.  This implies that there is some $k\in\Sigma_{d,h}\cap\Sigma_{d,h'}$.  Chasing the definitions, this implies that $h-\bar h = h' - \bar h' ($mod $ n)$.  In particular, if we define the sets
\begin{equation}\label{defd}
D_d^i:=\{j\in\{1,\dots,n-1\}:(i-1)c\leq j<ic\},
\end{equation}
then $h$ and $h'$ belong to the same $D_d^i$ and it follows that $\Sigma_{d,h}=\Sigma_{d,h'}$.  But for a fixed $i$, each element in $D_d^i$ has different reduction mod $c$.  Since $k_{\bar h}=k_{\bar h'}\in\Sigma_{d,h}=\Sigma_{d,h'}$, then we must have $\bar h= \bar h'$ implying that $h=h'$.
\end{proof}

We saw in the proof of Lemma \ref{lem:bij} that $h,h'\in D_d^i$ if and only if $\Sigma_{d,h}=\Sigma_{d,h'}$.  For this reason, we adopt the notation $\Sigma_d^i$.

We are now ready to cut down the rows in the matrix $\tilde\Phi$.  Using the above ordering convention, for any $\nu$ with $\nu^{tw}=\nu$ and $\underline\nu=\tau$, we can write
\[
\nu=((\tau_1,h_1),\dots,(\tau_l,h_l).
\]
We define $\hat\Phi$ to be the matrix obtained from $\tilde\Phi$ by restricting the row index to the set
\[
\{(\mu,k):m_1=0, k=k_{\tau_1}(h_1), \tilde\mu=-g_k(\tilde\nu) \text{ for some } \nu=\nu^{tw}\}.
\]
The fact that $\hat\Phi$ is square follows from Lemma \ref{lem:bij}.

\subsection{Invertibility of the Matrix}\label{subsec:invert}

To prove that $\hat\Phi$ is invertible over $\C((x))$, we proceed in two steps.  We first define certain blocks $\hat\Phi_h^i$ in $\hat\Phi$ with the following properties:
\begin{enumerate}
\item Each $\hat\Phi_h^i$ is invertible over $\C((x))$.
\item Each row and column of $\hat\Phi$ intersects exactly one $\hat\Phi_h^i$.
\item If $f(x)$ is an entry in some $\hat\Phi_h^i$ and $g(x)$ is an entry of $\hat\Phi$ in the same column, then $\ord_xf(x)\leq\ord_xg(x)$.
\end{enumerate}
If the inequality in (3) were strict, we would be done because the least degree term of the determinant of $\hat\Phi$ would be a signed product of the least degree terms in the determinants of the $\hat\Phi_h^i$ (by (2)) which are nonzero (by (1)).  However, the inequality  is not always strict as we will see below.  The second step is to use elementary matrix operations to take care of terms where the inequality is not strict.

We now define the blocks $\hat\Phi_h^i$.  For $h\in\{1,...,n-1\}^{l-1}$ and for $1\leq i\leq \frac{n}{c_1}$ define
\begin{eqnarray*}
B^i_h&=&\{ h_1\in D_{\tau_1}^i, t(\tilde\nu)=h \}\\
C^i_h&=&\{ k\in\Sigma_{\tau_1}^i, t(g_k(-\tilde\mu))=h \}.
\end{eqnarray*}
Then we define the sub-matrix $\hat\Phi^i_h$ by intersecting the indexing sets of $\hat\Phi$ with $B^i_h$ and $C^i_h$.

\begin{remark}\label{clarify2}
The above definitions might seem a bit obscure, a priori, but the motivation is simple.  From Remark \ref{clarify}, we know that the wreath Hurwitz numbers encoded by $\hat\Phi$ are rather simple.  In particular, the $\Z_n$ quotient of the cover is a disjoint union of $\proj^1$s and the only allowable monodromy over $\C^*\subset\proj^1$ are $x_1$ points.  For a fixed $\nu\in B^i_h$, the pairs $(\mu,k)\in C^i_h$ were chosen to be exactly those pairs such that there exists a wreath cover with the following three properties:
\begin{enumerate}
\item The $\Z_n$ monodromy over $0$ and $\infty$ for the $i$th $\proj^1$ is identified with $-h_i+\tau_ik$ and $m_i$, respectively,
\item The $\Z_n$ monodromy over the first $\C^*\subset\proj^1$ has the minimal possible number of $x_1$ points as we vary over all choices $(\mu,k)$ (this minimal number is $\bar h_1$), and
\item The $\Z_n$ monodromy over the other $\C^*$s is trivial.
\end{enumerate}
If we vary $\nu\in B^i_h$, the set of $(\mu,k)$ with these properties remains constant and they define the matrix $\hat\Phi^i_h$.
\end{remark}

\begin{remark}\label{clarify3}
That each column and each row of $\hat\Phi$ intersects exactly one $\hat\Phi_h^i$ follows from the fact that $D_{\tau_1}^i$ and $\Sigma_{\tau_1}^i$ both partition the set $\{1,\dots,n-1\}$.
\end{remark}

\begin{lemma}\label{lemma:det}
Let $\underline{\hat\Phi^i_h}$ denote the matrix of leading terms in $\hat\Phi^i_h$.  Then $\underline{\hat\Phi^i_h}$ is invertible.  In particular, $\hat\Phi^i_h$ is invertible over $\C((x))$.
\end{lemma}

\begin{proof}
By Remark \ref{clarify2}, the lowest degree term of the $((\mu,k),\nu)$ entry of $\hat\Phi_h^i$ has coefficient
\begin{equation}\label{wrh}
z_{\nu}\frac{(\xi_{2n}^{-1}\frac{k}{n})^{\bar h_1}}{\bar h_1!}H_{\mu,g_k(\nu)}^{2l,\gamma(h_1)\bullet}
\end{equation}
where $\bar{h}_1$ is \emph{independent} of $(\mu,k)\in C^i_h$ and $\gamma(h_1)$ is a $\bar h_1$-tuple of $\xi$'s.  The wreath Hurwitz numbers appearing in \eqref{wrh} are easy to compute, explicitly we have
\begin{align*}
z_{\nu}&\frac{(\xi_{2n}^{-1}\frac{k}{n})^{\bar h_1}}{\bar h_1!}H_{\mu,g_k(\nu)}^{2l,\gamma(h_1)\bullet}\\
&\hspace{1cm}=z_{\nu}\frac{(\xi_{2n}^{-1}\frac{k}{n})^{\bar h_1}}{\bar h_1!}\frac{1}{
|\aut(\mu)|}\nu_1^{\bar h_1}\prod_{i=1}^{l}\frac{1}{n\tau_i}\\
&\hspace{1cm}=\frac{|\aut(\nu)|}{|\aut(\mu)|}\frac{(\xi_{2n}^{-1}\frac{k}{n}\tau_1)^{\bar h_1}}{\bar h_1!}
\end{align*}

Therefore $\det\left(\underline{\hat\Phi^i_h}\right)$ is equal to
\[
\left(\prod_{(\mu,k)\in C^i_h}\frac{1}{|\aut(\mu)|}\right)\left(\prod_{\nu\in B^i_h}\frac{|\aut(\nu)|(\xi_{2n}^{-1}\nu_1x)^{\bar h_1}}{\bar h_1!}\right)
\det\left(\left(\frac{k}{n}\right)^{\bar h_1}\right)_{(\mu,k)\in C^i_h,\nu\in B^i_h}
\]
This is nonzero because $\det\left(\left(\frac{k}{n}\right)^{\bar h_1}\right)$ is the determinant of a Vandermonde matrix with different $k$ in different rows.
\end{proof}

Theorem \ref{thm:invert} now follows from the next result.

\begin{lemma}
$\hat\Phi$ is invertible over $\C((x))$.
\end{lemma}

\begin{proof}
For any fixed column $\alpha_{\nu}$ of $\hat\Phi$, there is a unique sub-matrix $\hat\Phi^i_h$ that intersects with this column. The degrees of the entries that lie in the intersection of $\alpha_\nu$ and $\hat\Phi^i_h$ are $\bar h_1$.  By the ordering convention introduced above, the degrees of the other entries of $\alpha_\nu$ are greater or equal to $\bar h_1$ (note that $m_1$ is always trivial). The equality holds for an entry in the row indexed by $(\mu,k)\notin C_h^i$ only if the following conditions are satisfied:
\begin{enumerate}
\item There exists a $j>1$ such that $\tau_j=\tau_1$, $\bar{h}_1=\bar{h}_j$, and $h_1<h_j$,
\item $-h_j+\tau_jk=-\bar{h}_j ($mod $n)$, and
\item $g_k(-\tilde{\mu})=\hat{\nu}$ where $\hat\nu=\nu\setminus\{ (\tau_j,h_j) \}$.
\end{enumerate}

If these conditions are met for some $(\mu,k)\notin C_h^i$, then there is a unique sub-matrix $\hat\Phi^{i'}_{h'}$ that intersects this row.  By definition, $h'=t(\hat\nu)$ and $i'$ is determined by the property $k\in\Sigma_{\tau_1}^{i'}$. It is not hard to see that every other entry that lies in the intersection of $\alpha_\nu$ and a row of $\hat\Phi^{i'}_{h'}$ also has minimal degree $\bar h_1$.
\[
\left(
\begin{array}{ccccccc}
&&&*&\cdots&\cdots&\cdots\\
&\hat\Phi^i_h&&*&\cdots&\cdots&\cdots\\
&&&*&\cdots&\cdots&\cdots\\
*&*&*&\ddots&*&*&*\\
*&*&x^{\bar h_1}&*&&&\\
\vdots&\vdots&\vdots&*&&\hat\Phi^{i'}_{h'}&\\
*&*&x^{\bar h_1}&*&&&
\end{array}\right)
\]

For every column $\alpha_{\nu'}$ that intersects $\hat\Phi^{i'}_{h'}$, we know $\widetilde{\nu'}=\hat\nu$.  In particular, $\bar{h}_2'=\bar{h}_1$ implying that $\bar{h}_1'\leq \bar{h}_1$ by the ordering convention.  If $\bar h_1'=\bar h_1$, then $\bar h_1'=\bar h_j$ (by (1)) and $h_1',h_j\in D_{\tau_1}^{i'}$ (the latter inclusion follows from (2)).  This would imply that $h_1'=h_j'$, i.e. $\nu=\nu'$ -- a contradiction.  Therefore we conclude that $\bar{h}_1'< \bar{h}_1=\bar{h}_2'$. In other words, condition (1) can never be satisfied by $\nu'$.  In particular, the degrees of the entries in $\alpha_{\nu'}$ which are not contained in $\hat\Phi^{i'}_{h'}$ are strictly greater than $\bar h_1'$.

By Lemma \ref{lemma:det}, we can transform the matrix $\hat\Phi^{i'}_{h'}$ to a matrix $\Psi^{i'}_{h'}$ such that $\Psi^{i'}_{h'}|_{x=0}$ is the identity matrix.  More specifically, we first multiply each column by $x^{-\bar h_1'}$ where $\nu'$ is the index of the column, then we apply elementary column operations (over $\C$) to reduce the matrix of (constant) leading terms to the identity.  Extending these column operations to the columns of $\hat\Phi$, we can replace the sub-matrix $\hat\Phi^{i'}_{h'}$ by $\Psi^{i'}_{h'}$ in such a way that the following two properties are satisfied:
\begin{enumerate}[(a)]
\item For each column intersecting $\Psi^{i'}_{h'}$, the entries which do not lie in $\Psi^{i'}_{h'}$ have vanishing constant terms, and
\item The transformed matrix is invertible over $\C((x))$ if and only if the original matrix is invertible over $\C((x))$.
\end{enumerate}

\[
\left(
\begin{array}{ccccccc}
&&&*&&&\\
&\hat\Phi^i_h&&*&&O(x)&\\
&&&*&&&\\
*&*&*&\ddots&*&*&*\\
*&*&x^{\bar h_1}&*&&&\\
\vdots&\vdots&\vdots&*&&\Psi^{i'}_{h'}=I+O(x)&\\
*&*&x^{\bar h_1}&*&&&
\end{array}\right)
\]

We can now use the columns intersecting $\Psi^{i'}_{h'}$ to cancel the degree $\bar h_1$ terms of the entries that lie in the intersection of $\alpha_\nu$ and rows of $\Psi^{i'}_{h'}$.  By property (a), this does not affect the degree $\bar h_1$ terms in the entries of $\alpha_\nu$ which do not lie in rows which intersect $\Psi^{i'}_{h'}$.  In particular, the smallest degree term in $\det\left(\hat\Phi_{h}^i \right)$ is not affected.  We can repeat this process until the least degree terms in each column are contained in the sub-matrix $\hat\Phi_{h}^i$ (or $\Psi_{h}^i$ if it has been transformed).  Call the resulting matrix $\Psi$.  Then the least degree term of $\det(\Psi)$ is the product of least degree terms of determinants of matrices of the form $\hat\Phi_{h}^i$ or $\Psi_{h}^i$, all of which are nonzero.  Therefore $\Psi$ is invertible over $\C((x))$.  By property (b), $\hat{\Phi}$ is invertible over $\C((x))$.\end{proof}

\section{Combinatorics}\label{sec:combinatorics}

In this section, we investigate the framed Donaldson-Thomas vertex $\tilde{P}_\lambda(a)$ and relate it to loop Schur functions.

\subsection{Loop Schur Functions}

For a positive integer $n$ and partition $\rho$, the \textit{colored} Young diagram $(\rho,n)$ is obtained by coloring the boxes of the Young diagram by their \textit{content} modulo $n$.  In other words if $\square$ is in the $i$th row and the $j$th column, we color it with $c(\square):=j-i \mod n$.  For example, if $\rho=(4,3,3,1)$ and $n=3$, the colored Young diagram is given by
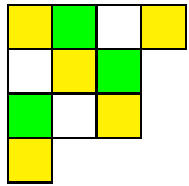
\begin{center}
\begin{ytableau}
*(yellow) & *(green) & *(white) & *(yellow)\\
*(white) & *(yellow) & *(green)\\
*(green) & *(white) & *(yellow)\\
*(yellow)
\end{ytableau}
\end{center}
with
\[
0\leftrightarrow\begin{ytableau}*(yellow)\end{ytableau}, \hspace{1cm} 1\leftrightarrow\begin{ytableau}*(green)\end{ytableau}, \hspace{.5cm}\text{and}\hspace{.5cm} 2\leftrightarrow\begin{ytableau}*(white)\end{ytableau}
\]

We let $\rho[i]$ denote the collection of boxes with color $i$.  A \textit{semi-standard Young tableau} (SSYT) of $\rho$ is a numbering of the boxes so that numbers are weakly increasing left to right and strictly increasing top to bottom.   For each SSYT $T$ and $\square\in\rho$, we define the \textit{weight} $w(\square,T)$ to be the number appearing in that box.  To each $\rho$, $n$, and $T\in SSYT(\rho,n)$, we associate a monomial $q^T$ in $n$ infinite sets of variables $\{q_{i,j}|i\in\Z_n,j\in\N\}$:
\[
q^T:=\prod_{i=0}^{n-1}\prod_{\square\in\rho[i]}q_{i,w(\square,T)}.
\]
For example, to the SSYT
\[
T=\begin{ytableau}
*(yellow) 1& *(green) 1& *(white) 2& *(yellow) 4\\
*(white) 2& *(yellow) 3& *(green) 3\\
*(green) 4& *(white) 4& *(yellow) 6\\
*(yellow) 7
\end{ytableau}
\]
we associate the monomial
\[
q^T=q_{0,1}q_{0,3}q_{0,4}q_{0,6}q_{0,7}q_{1,1}q_{1,3}q_{1,4}q_{2,2}^2q_{2,4}.
\]
\begin{definition}
The \textit{loop Schur function} associated to $(\rho,n)$ is defined by
\[
s_{\rho}[n]:=\sum_{T\in SSYT(\rho,n)}q^T.
\]
\end{definition}

In the current setting, we are only concerned with the case where $\rho=\bar\lambda$ arises from an $n$-tuple of partitions $\lambda$ via $n$-quotients (c.f. Section \ref{sec:wrfock}).  This is equivalent to the following condition.

\begin{definition}\label{def:balanced}
We call the colored Young diagram $\rho$ \textit{balanced} if $|\rho[i]|=|\rho[j]|$ for all $i,j$.
\end{definition}

Denote by $S_{\lambda}$ the function in $n$ variables obtained by making the substitution $q_{i,j}=q_i^j$ in $s_{\bar\lambda}[n]$.  The following result is due to E\u{g}ecio\u{g}lu and Remmel \cite{er:cpgisf} and also appears in a more recent paper of Nakada \cite{n:hfgt}.

\begin{lemma}[\cite{er:cpgisf,n:hfgt}]\label{hookcontent}
\[
S_\lambda=\frac{\prod_i q_i^{n_i(\bar\lambda)}}{\prod_{\square\in\bar\lambda}\left(1-\prod_i q_i^{h_i(\square)}\right)}.
\]
\end{lemma}

As a consequence, we have the following identity:

\begin{corollary}\label{sp}
\[
\tilde{P}_\lambda(0)=\frac{\cha_{\bar\lambda}(n^d)}{\dim(\lambda)}q^\frac{d}{2}(-1)^dS_{\lambda}.
\]
\end{corollary}

We also recall the definition of the series $s_\rho^k[n]$ from \cite{r:cmn}.  For $0\leq k<n$, define the shifted weight
\[
w^k(\square,T):=w(\square,T)+\frac{k\cdot c(\square)}{n}
\]
and the corresponding monomial
\[
q^{T,k}:=\prod_{i=0}^{n-1}\prod_{\square\in\rho[i]}q_{i,w^k(\square,T)}
\]
where the second index belongs to $\frac{1}{n}\Z$.

\begin{definition} The \textit{$k$-shifted Schur function}\footnote{The $k$-shifted Schur functions here should not be confused with the shifted Schur functions defined by Okounkov and Olshanski \cite{oo:ssf}.  We shift the index of the variables wheras they shift the variables themselves.  Moreover, they sum over \textit{reverse} tableaux.} associated to $(\rho,n)$ is
\[
s_\rho^k[n]:=\sum_{T\in SSYT(\rho,n)}q^{T,k}.
\]
\end{definition}

We denote by $S_\lambda^k$ the series in $n$ variables obtained from $s_{\bar\lambda}^k[n]$ by specializing $q_{i,j}=q_i^j$.

\begin{remark}
Notice the specialization $s_\rho^0[n]=s_\rho[n]$, and hence similarly with $S$.
\end{remark}

Since $S_\lambda^k$ differs from $S_\lambda$ only by a monomial factor, we have the following natural generalization of Corollary \ref{sp}.

\begin{lemma}\label{sp2}
\[
\tilde{P}_\lambda(0)=\frac{\cha_{\bar\lambda}(n^d)}{\dim(\lambda)}q^\frac{d}{2}(-1)^dS_{\lambda}^k\left(\prod_{(i,j)\in\bar\lambda}q_{j-i}^{j-i}\right)^{\frac{-k}{n}}.
\]
\end{lemma}

\subsection{Combinatorial Identities}

Before stating the necessary combinatorial identities, we provide the following definition.

\begin{definition}
The set theoretic difference $\rho\setminus\tau$ of two Young diagrams $\rho\supset\tau$ is a \textit{border strip of $\tau$} if it is connected and does not contain any $2\times 2$ block.  The \textit{length} of a border strip is the number of boxes it contains.  The \textit{height} (denoted $ht$) of a border strip is the number of rows it occupies, minus $1$.
\end{definition}

The following are specializations of results proved by the first author \cite{r:cmn}.
\begin{theorem}[\cite{r:cmn}, Theorem 1]\label{thm:comb1}
\[
\frac{1}{1-(q_0...q_{n-1})^l}S_\lambda=\sum (-1)^{ht(\bar\sigma\setminus\bar\lambda)}S_\sigma
\]
where the sum is over all ways of adding a length $ln$ border strip to $\bar\lambda$.
\end{theorem}

\begin{theorem}[\cite{r:cmn}, Theorem 2]\label{thm:comb2}
For a fixed $\bar\lambda$ and $k\neq 0$,
\[
\sum (-1)^{ht(\bar\sigma\setminus\bar\lambda)}S_\sigma^k=0
\]
where the sum is over all ways of adding a length $ln$ border strip to $\bar\lambda$.
\end{theorem}

\section{Representation Theory}\label{sec:reptheory}

In this section we investigate certain characters of the generalized symmetric group which arose in Section \ref{sec:wreath}.  Our main tool is the wreath Fock space.  We begin by recalling the basic definitions and results concerning the usual Fock space.

\subsection{The Infinite Wedge}\label{sec:infwedge}

The infinite wedge provides a convenient setting for studying the representation theory of the symmetric group in terms of combinatorial manipulations of partitions and Maya diagrams.  For a more thorough treatment of the infinite wedge and some of its applications in Gromov-Witten theory, see for example the work of Okounkov and Pandharipande \cite{op:hiaiotu,op:egwtop} or for an application in double Hurwitz numbers, see Johnson \cite{j:dhviw}.

Let $V$ be the infinite vector space with spanning set indexed by half integers:
\[
V:=\bigoplus_{i\in\Z}\left\langle\underline{i+\frac{1}{2}}\right\rangle_{\C}.
\]
\begin{definition}
The \textit{infinite wedge} $\bigwedge^{\frac{\infty}{2}}V$ is the vector space
\[
\bigwedge^{\frac{\infty}{2}}V:=\bigoplus_{(i_k)}\left\langle\underline{i_1}\wedge\underline{i_2}\wedge...\right\rangle
\]
where $(i_k)$ is a decreasing sequence of half integers such that
\[
i_k+k-\frac{1}{2}=c
\]
for some constant $c$ and $k\gg0$.  We call $c$ the \textit{charge} of the vector.
\end{definition}

We will only be concerned with the subvector space spanned by vectors of charge $0$.  We denote this space by $\bigwedge_0^{\frac{\infty}{2}}V$.

\subsubsection{Maya Diagrams}

The primary combinatorial tool for us will be Maya diagrams.  A Maya diagram is a collection of stones placed at the half integers such that the half integers without stones are bounded below and the half integers with stones are bounded above.  A Maya diagram has \textit{charge zero} if the number of stones at positive half integers is equal to the number of negative half integers without stones.

The basis vectors of $\bigwedge_0^{\frac{\infty}{2}}V$ can be identified with charge zero Maya diagrams canonically as follows.  Let $S=\{i_k\}$ where $(i_k)$ corresponds to a charge $0$ vector.  Then we obtain a charge zero Maya diagram by placing a stone in the $i$th place if and only if $i\in S$.

\subsubsection{Partitions}

The charge zero basis vectors can also be canonically identified with partitions.  If we let $\alpha$ be the increasing sequence of half integers in $S\cap\Q_{>0}$ and $\beta$ the increasing sequence of half integers in $-(S^c\cap\Q_{<0})$, then $(\alpha | \beta)$ is the modified Frobenius coordinate of a partition $\rho$.  In other words, representing $\rho$ as a Young diagram, $\alpha_i$ is the number of boxes (half-boxes included) in the $i$th row to the right of the main diagonal and $\beta_i$ is the number of boxes in the $i$th column below the main diagonal.

Equivalently, the partition $\rho=(\rho_1,\rho_2,...)$ is determined by writing the vector $v_S$ in the following form.
\[
v_S=\underline{\rho_1-1/2}\wedge\underline{\rho_2-3/2}\wedge.....
\]

To relate partitions to Maya diagrams, rotate the corresponding Young diagram counterclockwise by 135 and place $0$ directly below the vertex.  The stones in the Maya diagram lie directly below outer edges of the Young diagram which have slope $1$.  This correspondence is illustrated in Figure 1.

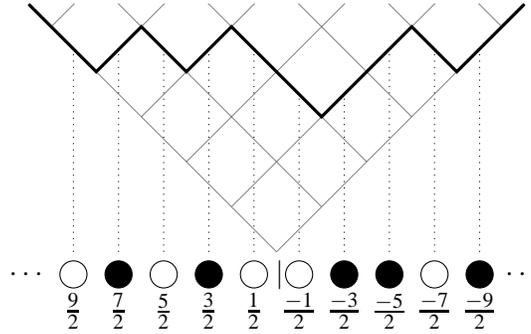
\begin{figure}\label{fig:bases}
\begin{center}
\begin{tikzpicture}
\begin{scope}[gray, very thin, scale=.6]
\clip (-5.5, 5.5) rectangle (5.5, -5);
\draw[rotate=45, scale=1.412] (0,0) grid (6,6);
\end{scope}
\begin{scope}[rotate=45, very thick, scale=.6*1.412]
\draw (0,5.5) -- (0, 4) -- (1,4) -- (1,3)-- (2,3) -- (2,1) -- (4,1) -- (4,0) -- (5.5,0);
\end{scope}
\begin{scope}[scale=.6, dotted]
\draw (-4.5,0) -- (-4.5, 4.5);
\draw (-3.5,0) -- (-3.5, 4.5);
\draw (-2.5,0) -- (-2.5, 4.5);
\draw (-1.5,0) -- (-1.5, 4.5);
\draw (-.5,0) -- (-.5, 4.5);
\draw (.5,0) -- (.5, 3.5);
\draw (1.5,0) -- (1.5, 3.5);
\draw (2.5,0) -- (2.5, 4.5);
\draw (3.5,0) -- (3.5, 4.5);
\draw (4.5,0) -- (4.5, 4.5);
\end{scope}
\begin{scope}[scale=.6, yshift=-.5cm]
\draw (-5.5,0) node{$\cdots$};
\draw (-4.5,0) circle (.3) node[below=3pt]{$\frac{9}{2}$};
\filldraw (-3.5,0) circle (.3) node[below=3pt]{$\frac{7}{2}$};
\draw (-2.5,0) circle (.3) node[below=3pt]{$\frac{5}{2}$};
\filldraw (-1.5,0) circle (.3) node[below=3pt]{$\frac{3}{2}$};
\draw (-.5,0) circle (.3) node[below=3pt]{$\frac{1}{2}$} node[right=3.5pt]{$|$};
\draw (.5,0) circle (.3) node[below=3pt]{$\frac{-1}{2}$};
\filldraw (1.5,0) circle (.3) node[below=3pt]{$\frac{-3}{2}$};
\filldraw (2.5,0) circle (.3) node[below=3pt]{$\frac{-5}{2}$};
\draw (3.5,0) circle (.3) node[below=3pt]{$\frac{-7}{2}$};
\filldraw (4.5,0) circle (.3) node[below=3pt]{$\frac{-9}{2}$};
\draw (5.5,0) node{$\cdots$};
\end{scope}
\end{tikzpicture}
\end{center}
\caption{Correspondence between the different combinatorial bases of $\bigwedge_0^{\frac{\infty}{2}}V$.}
\end{figure}

\subsubsection{One Basis}

With the above correspondences, we will think of $\bigwedge_0^{\frac{\infty}{2}}V$ simultaneously as the vector space spanned by
\begin{itemize}
\item Sequences $S$ of the half integers with charge $0$,
\item Maya diagrams with charge $0$, or
\item Partitions.
\end{itemize}
For simplicity, we will denote the basis elements by $v_\rho$ keeping in mind that the partition $\rho$ corresponds canonically to a Maya diagram $m_\rho$ and a set of half integers $S_\rho$.  We denote by $v_\emptyset$ the \textit{vacuum vector} which is the vector corresponding to the trivial partition.

\subsubsection{Operators}

In order to relate the infinite wedge to the representation theory of $S_d$, we define several operators on $\bigwedge_0^{\frac{\infty}{2}}V$ via their action on basis elements $v_\rho$.

For any half integer $k$ and basis element $v_\rho$, the operator $E_{k,k}$ acts on $v_\rho$ as follows:
\[
E_{k,k}v_\rho=\begin{cases}
v_\rho &k>0,k\in S_\rho\\
-v_\rho &k<0,k\notin S_\rho\\
0 &\text{else.}
\end{cases}
\]

For $k$ a positive integer, the creation operator $\alpha_{-k}$ acts on $v_\rho$ as follows:
\[
\alpha_{-k}v_\rho=\sum_{\tau}(-1)^{ht(\tau\setminus\rho)}v_{\sigma}
\]
where the sum is over all ways of adding length $k$ border strips to $\rho$.  In terms of Maya diagrams, the sum is over all ways of moving a stone $k$ places to the left and the sign corresponds to the number of stones jumped during such a move.

Recall that each partition $\rho$ corresponds to an irreducible representation of $S_d$ with character $\cha_\rho$.  Given a partition $\tau=(d_1,...,d_l)$ corresponding to a conjugacy class in $S_d$, we define the operator
\[
\alpha_{-\tau}:=\prod_{i=1}^l\alpha_{-d_i}
\]
The following identity follows from the Murnaghan-Nakayama formula.
\begin{equation}\label{eqn:create}
\alpha_{-\tau}v_\emptyset=\sum_\rho\cha_\rho(\tau)v_\rho.
\end{equation}

We also define the operator
\[
\mathcal{F}_T:=\sum_k\frac{k^2}{2}E_{k,k}.
\]
If $T$ is the conjugacy class of transpositions and $f_T(\lambda):=\frac{|T|\cha_\lambda(T)}{\dim(\lambda)}$, then each $v_\lambda$ is an eigenvector of $\mathcal{F}_T$ with eigenvalue $f_T(\lambda)$:
\begin{equation}\label{eqn:central}
\mathcal{F}_T\cdot v_\lambda=f_T(\lambda)v_\lambda.
\end{equation}

\subsection{Wreath Fock Space}\label{sec:wrfock}

The wreath product generalization of the Fock space gives a combinatorial tool for manipulating the representation theory of the groups $G\wr S_d$.  These spaces and their corresponding operators have been developed in e.g. \cite{fw:vawpc,qw:hsop,j:egwtods}.  We merely focus on the cyclic case which is all we require.  To that end, the wreath Fock space can be defined as
\[
\mathcal{Z}_n:=\bigotimes_{\{0,...,n-1\}}\bigwedge\vphantom{1cm}_{\vspace{-.1cm}0\vspace{.1cm}}^{\frac{\infty}{2}}V.
\]
Basis vectors correspond to $n$-tuples of partitions $\lambda=(\lambda_0,...,\lambda_{n-1})$ or, equivalently, $n$-tuples of Maya diagrams.

In the wreath Fock space, there is an additional way by which we will distinguish a basis element.  Given an $n$-tuple of Maya diagrams, we can interlace them to get a single Maya diagram by sending a stone in the $k$th place of the $i$th Maya diagram to position $n\left(k-\frac{1}{2}\right)+(i+\frac{1}{2})$ in the new Maya diagram.  An example of this identification is shown in Figure 2.  This new Maya diagram corresponds to a partition of $nd$ which we denote $\bar\lambda$.  Reversing this process is usually referred to as an \textit{$n$-quotient}.  It is well known that taking $n$-quotients gives a bijection between balanced Young diagrams $\bar\lambda$ (c.f. Definition \ref{def:balanced}) and $n$-tuples of partitions $\lambda$.

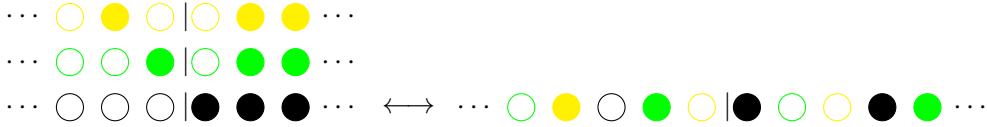
\begin{figure}\label{fig:maya}
\begin{tikzpicture}
\begin{scope}[scale=.6, yshift=-.5cm]
\draw (-1.5,0) node{$\cdots$};
\draw (-.5,0) circle (.3);
\draw (.5,0) circle (.3);
\draw (1.5,0) circle (.3) node[right=3.5pt]{\color{black}$|$};
\filldraw (2.5,0) circle (.3);
\filldraw (3.5,0) circle (.3);
\filldraw (4.5,0) circle (.3);
\draw (5.5,0) node{$\cdots$};
\draw (-1.5,1) node{$\cdots$};
\draw[green] (-.5,1) circle (.3);
\draw[green] (.5,1) circle (.3);
\filldraw[green] (1.5,1) circle (.3) node[right=3.5pt]{\color{black}$|$};
\draw (2.5,1)[green] circle (.3);
\filldraw[green] (3.5,1) circle (.3);
\filldraw[green] (4.5,1) circle (.3);
\draw (5.5,1) node{$\cdots$};
\draw (-1.5,2) node{$\cdots$};
\draw[yellow] (-.5,2) circle (.3);
\filldraw[yellow] (.5,2) circle (.3);
\draw[yellow] (1.5,2) circle (.3) node[right=3.5pt]{\color{black}$|$};
\draw[yellow] (2.5,2) circle (.3);
\filldraw[yellow] (3.5,2) circle (.3);
\filldraw[yellow] (4.5,2) circle (.3);
\draw (5.5,2) node{$\cdots$};
\draw (7,0) node{$\longleftrightarrow$};
\draw (8.5,0) node{$\cdots$};
\draw[green] (9.5,0) circle (.3);
\filldraw[yellow] (10.5,0) circle (.3);
\draw (11.5,0) circle (.3);
\filldraw[green] (12.5,0) circle (.3);
\draw[yellow] (13.5,0) circle (.3) node[right=3.5pt]{\color{black}$|$};
\filldraw (14.5,0) circle (.3);
\draw[green] (15.5,0) circle (.3);
\draw[yellow] (16.5,0) circle (.3);
\filldraw (17.5,0) circle (.3);
\filldraw[green] (18.5,0) circle (.3);
\draw (19.5,0) node{$\cdots$};
\end{scope}
\end{tikzpicture}
\caption{A example of a $3$-quotient.}
\end{figure}

For any operator $M$ on $\bigwedge_0^{\frac{\infty}{2}}V$ and any integer $0\leq k\leq n-1$, we define the operator $M^k$ to act on $\mathcal{Z}_n$ by acting as $M$ on the $k$th factor and trivially on the other factors.

Given $\lambda$, we can canonically identify it with an irreducible representation of $\Z_n\wr S_d$ with character $\cha_\lambda$.  Similarly, given an $n$-tuple of partitions $\mu=(\mu^0,...,\mu^{n-1})$ with $\mu^k=(d_1^k,...,d_{l_k}^k)$, we can be canonically identify it with a conjugacy class.  We have the following important generalizations of (\ref{eqn:create}) and (\ref{eqn:central}):
\begin{equation}\label{eqn:create2}
\prod_{k=0}^{n-1} \prod_{i=0}^{l_k}\left(\sum_{j=0}^{n-1} \xi^{-kj}\alpha_{-d_i^k}^j\right) v_\emptyset=\sum_\lambda\cha_\lambda(\mu)v_\lambda.
\end{equation}
and
\begin{equation}\label{eqn:central2}
\left(n\sum_{i=0}^{n-1}\mathcal{F}_T^i\right)\cdot v_\lambda=f_T(\lambda)v_\lambda.
\end{equation}

\subsection{Central Characters}

We now use the combinatorics of colored partitions and Maya diagrams to study the central characters $f_i(\lambda)$ and $f_T(\lambda)$ which arose in Section \ref{sec:wreath}.

\begin{lemma}\label{lem:central}
Let $\lambda=(\lambda_0,...,\lambda_{n-1})$ with $|\lambda_i|=d_i$.  Then
\begin{enumerate}[(i)]
\item \[f_i(\lambda)=\sum_j\xi_n^{-ij}d_j\]\label{ident:i} and
\item \[f_T(\lambda)=\sum_{(i,j)\in\bar\lambda[0]}j-i\]\label{ident:ii}
\end{enumerate}
\end{lemma}
\begin{proof}
To prove identity \eqref{ident:i}, recall that
\[
f_i(\lambda)=\frac{d\cha_\lambda(\{\xi^i,1^{d-1}\})}{\dim(\lambda)}=\frac{d\cha_\lambda(\{\xi^i,1^{d-1}\})}{\cha_\lambda(\{1^d\})}
\]
where the exponent of $1$ in the multiset denotes repetition.  For $\mu=\{1^d\}$, the coefficient of $v_\lambda$ in (\ref{eqn:create2}) can be interpreted as the number of ways to build the $n$-tuple of Young diagrams $\lambda=(\lambda_0,...,\lambda_{n-1})$ one box at a time.  Equivalently, this can be interpreted as the number of standard Young tableaux of $\lambda$, i.e. the number of ways to fill the boxes of the $\lambda_i$ with the numbers $1,...,d$ with the property that numbers always increase along rows and down columns.  This is easily computed:
\begin{equation}\label{cc1}
\cha_\lambda(\{1^d\})={d \choose d_0,...,d_{n-1}}\prod\dim(\lambda_i)
\end{equation}
where we use the fact that $\dim(\lambda_i)$ is the number of standard tableaux of $\lambda_i$.

On the other hand, for $\mu=\{\xi^i,1^{d-1}\}$, the coefficient of $v_\lambda$ in (\ref{eqn:create2}) can be interpreted as a weighted count of ways to build $\lambda$ one box at a time, where the weight is $\xi^{-ij}$ if the first box is a part of $\lambda_j$.  This is also easily computed:
\begin{equation}\label{cc2}
\cha_\lambda(\{\xi^i,1^{d-1}\})=\sum_{j=0}^{n-1}\xi_n^{-ij}{d-1 \choose d_0,...,d_j-1,...,d_{n-1}}\prod\dim(\lambda_i).
\end{equation}
Identity \eqref{ident:i} follows by dividing \eqref{cc2} by \eqref{cc1} and multiplying by $d$.

To prove identity \eqref{ident:ii}, begin by writing $\bar\lambda=(\alpha|\beta)$ in modified Frobenius notation (c.f. Section \ref{sec:infwedge}).  Then the number of boxes in $\bar\lambda[0]$ to the right (below) the $i$th diagonal element is given by $\left\lfloor \frac{\alpha_i}{n} \right\rfloor$ $\left(\left\lfloor \frac{\beta_i}{n} \right\rfloor\right)$.  If we compute the sum in \eqref{ident:ii} over these $\left\lfloor \frac{\alpha_i}{n} \right\rfloor$ $\left(\left\lfloor \frac{\beta_i}{n} \right\rfloor\right)$ terms, we get a contribution of
\[
n+2n+...+n\left\lfloor \frac{\alpha_i}{n} \right\rfloor \hspace{1cm}\left(-n-2n-...-n\left\lfloor \frac{\beta_i}{n} \right\rfloor  \right).
\]

Therefore, the right side of the \eqref{ident:ii} can be written as
\begin{equation}\label{eqn:rhs}
\sum_{(i,j)\in\bar\lambda[0]}j-i=n\sum_{i=1}^{\infty}\left(\frac{\left\lfloor \frac{\alpha_i}{n} \right\rfloor^2+\left\lfloor \frac{\alpha_i}{n} \right\rfloor}{2}- \frac{\left\lfloor \frac{\beta_i}{n} \right\rfloor^2+\left\lfloor \frac{\beta_i}{n} \right\rfloor}{2} \right).
\end{equation}

To compute the left side of \eqref{ident:ii}, we consider equation (\ref{eqn:central2}).  Via the $n$-quotient correspondence described above, we can interpret $v_\lambda$ as a vector $v_{\bar\lambda}\in\bigwedge_0^{\frac{\infty}{2}}$.  Under this correspondence, the operator $n\sum_{i=0}^{n-1}\mathcal{F}_T^i$ becomes
\[
n\sum_k\frac{1}{2}\left(\left\lfloor\frac{k}{n}\right\rfloor+\frac{1}{2}\right)^2E_{kk}.
\]
Each summand acts simply by multiplying $v_{\bar\lambda}$ by an appropriate scalar.  This scalar is zero unless $k=\alpha_i>0$ or $k=-\beta_i<0$ for some $i$.  In these cases, the scalar is
\[
n\frac{1}{2}\left(\left\lfloor\frac{\alpha_i}{n}\right\rfloor+\frac{1}{2}\right)^2
\]
and
\[
-n\frac{1}{2}\left(\left\lfloor\frac{\beta_i}{n}\right\rfloor+\frac{1}{2}\right)^2.
\]
We obtain (\ref{eqn:rhs}) by summing over all such $i$.
\end{proof}

\begin{lemma}\label{lem:frame}
After the change of variables prescribed by Theorem \ref{thm1},
\begin{equation}\label{eqn:toshow}
\left(\prod_{(i,j)\in\bar\lambda} q_{j-i}^{j-i}\right)^{1/n}=\left(-\xi_{2n}\right)^{-d} \left(\xi_n^{-\sum kd_k}\right)\re^{\frac{1}{n}\left(\sqrt{-1}f_T(\lambda)u+\sum \xi_{2n}^{-k}f_k(\lambda)x_k\right)}
\end{equation}
\end{lemma}

\begin{proof}
If $\lambda=(\lambda_0,...,\lambda_{n-1})$ with $|\lambda_k|=d_k$, then in terms of Maya diagrams we can interpret the $d_i$ as follows: $d_k$ is the number of moves it takes to build the Maya diagram of $\lambda_k$ from the empty Maya diagram by only moving stones one place at a time.  Moreover, each such move has the effect of adding a length $n$ border strip to $\bar\lambda$, the northeast-most box in the strip having color $k$.  The quantity $j-i$ decreases uniformly by $1$ as we move south and west along the strip so each such move contributes to $\prod_{(i,j)\in\bar\lambda} q_{j-i}^{j-i}$ a factor of
\begin{equation}\label{stripcontrib}
q_k^lq_{k-1}^{l-1}\cdots q_1^{l-k+1}q_0^{l-k}q_{n-1}^{l-k-1}\cdots q_{k+1}^{l-n+1}
\end{equation}
for some $l$.  In order to apply the change of variables, we need to collect the $q_0$'s into $q$'s.  Borrowing the necessary $q_i$'s from the other squares in the border strip, \eqref{stripcontrib} becomes
\[
q^{l-k}\left( q_k^kq_{k-1}^{k-1}\cdots q_1^{1}q_{n-1}^{-1}\cdots q_{k+1}^{k-n+1} \right).
\]
Combining these factors for all $k$, we find
\begin{equation}\label{eqn:factors}
\prod_{(i,j)\in\bar\lambda} q_{j-i}^{j-i}=q^M\prod_{k=0}^{n-1}\left( q_k^kq_{k-1}^{k-1}\cdots q_1^{1}q_{n-1}^{-1}\cdots q_{k+1}^{k-n+1} \right)^{d_k}
\end{equation}
where $M=\sum_{(i,j)\in\bar\lambda[0]}(j-i)$ is the total power of $q_0$ which we know is equal to $f_T(\lambda)$ from Lemma \ref{lem:central}.

It is left to investigate what happens to the factors in (\ref{eqn:factors}) after the change of variables.  Since $q\rightarrow \re^{\sqrt{-1}u}$ and $M=f_T(\lambda)$, then we see immediately that the $u$ factors on either side of (\ref{eqn:toshow}) agree.

We now compute the coefficient of $d_ix_j$ in the exponent of (\ref{eqn:factors}) after the change of variables.  To do this, we must compute the coefficient of $x_j$ in the factor $q_i^iq_{i-1}^{i-1} \cdots q_1^{1}q_{n-1}^{-1}\cdots q_{i+1}^{i-n+1}.$  Applying the change of variables, this coefficient is
\begin{equation}\label{eq:exp}
-\sum_{r=1}^i\frac{r\xi_n^{-jr}}{n}\left( \xi_{2n}^j-\xi_{2n}^{-j} \right)-\sum_{s=i+1}^{n-1}\frac{(s-n)\xi_n^{-js}}{n}\left( \xi_{2n}^j-\xi_{2n}^{-j} \right).
\end{equation}
Setting $y:=\xi_n^{-j}$, \eqref{eq:exp} can be written as
\begin{align}\label{eqn:factors2}
\nonumber\frac{-y^{-\frac{1}{2}}}{n}&\left(\sum_{r=1}^i \left(ry^r-ry^{r+1}\right)+\sum_{s=i+1}^{n-1}\left((s-n)y^s-(s-n)y^{s+1}\right) \right)\\
&\hspace{1cm}=\frac{-y^{-\frac{1}{2}}}{n}\left(-ny^{i+1}+ \sum_{r=1}^n y^r \right).
\end{align}

Using the fact that $\sum_{r=1}^n y^r=0$,  \eqref{eqn:factors2} is equal to $\xi_{2n}^{j(-2i-1)}$.  Therefore, the coefficient of $x_j$ is
\[
\xi_{2n}^{-j} \sum \xi_n^{-ij}d_i=\xi_{2n}^{-j} f_j(\lambda)
\]
where the equality follows from the first identity of Lemma \ref{lem:central}.

Finally, notice that the root of unity which factors out of the term \[\left( q_i^iq_{i-1}^{i-1}\cdots q_1^{1}q_{n-1}^{-1}\cdots q_{i+1}^{i-n+1} \right)^{1/n}\] after the change of variables is $-\xi_{2n}^{-1}\xi_n^{-i}$.  Putting all of this together proves the result.

\end{proof}

\subsection{Signs}\label{sec:signs}

If $\bar\sigma$ is obtained from $\bar\lambda$ by adding a length $kn$ border strip, then the Maya diagrams corresponding to $\sigma$ are obtained from those corresponding to $\lambda$ by moving a stone $k$ places in the $i$th Maya diagram.  Notice that $k$ and $i$ are both determined by $\bar\sigma$ and $\bar\lambda$.  For notational convenience, we make the following definition.

\begin{definition}
If $\bar\sigma$ is obtained from $\bar\lambda$ by adding a length $kn$ border strip, let $\beta(\sigma\setminus\lambda)$ denote the number of stones in the $i$th Maya diagram which are skipped over.
\end{definition}

Notice that $(-1)^{\beta(\sigma\setminus\lambda)}$ is the coefficient of $v_\sigma$ in $\alpha_{-k}^i(v_\lambda)$.

The next lemma allows us to deal with the sign $\frac{\cha_{\bar\lambda}(n^d)}{\dim(\lambda)}$ appearing in Theorem \ref{thm1}.

\begin{lemma}\label{signs1}
If $\bar\sigma$ is obtained from $\bar\lambda$ by adding a length $kn$ border strip, then
\[
\frac{\cha_{\bar\sigma}(n^{d+k})}{\dim(\sigma)}=(-1)^{\beta(\sigma\setminus\lambda)+ht(\bar\sigma\setminus\bar\lambda)-1}\frac{\cha_{\bar\lambda}(n^d)}{\dim(\lambda)}.
\]
\end{lemma}

\begin{proof}
By (\ref{eqn:create}), $\cha_{\bar\lambda}(n^d)$ is the weighted sum of ways to create the Maya diagram of $\bar\lambda$ from the vacuum diagram by moving stones $n$ places at a time; the weight is $\pm 1$ depending on whether the total number of stones jumped over is even or odd.  It is not hard to see that the weight of any such sequence is equal to the weight of any other.  Since $\dim(\lambda)$ is the total number of such sequences, we see that $\frac{\cha_{\bar\lambda}(n^d)}{\dim(\lambda)}$ is equal to the weight of any one of them.

Now suppose $\bar\sigma$ is obtained from $\bar\lambda$ by adding a length $kn$ border strip.  We can think of $\bar\sigma$ as being obtained from $\bar\lambda$ by moving a single stone $kn$ places to the left in the Maya diagram of $\bar\lambda$, $ht(\bar\sigma\setminus\bar\lambda)-1$ is the total number of stones jumped while $\beta(\sigma\setminus\lambda)$ counts the number of jumped stones which are $n,2n,3n,...$ positions to the left of where the stone sat originally.

On the other hand, the Maya diagram of $\bar\sigma$ can be obtained from that of $\bar\lambda$ by choosing a sequence of length $n$ jumps.  As above, $\frac{\cha_{\bar\sigma}(n^{d+k})}{\dim(\sigma)}=(-1)^*\frac{\cha_{\bar\lambda}(n^d)}{\dim(\lambda)}$ where $*$ is equal to the total number of stones jumped during the sequence of moves.  With the above interpretations for $ht(\bar\sigma\setminus\bar\lambda)-1$ and $\beta(\sigma\setminus\lambda)$, we see that the number of stones jumped in this process is $(ht(\bar\sigma\setminus\bar\lambda)-1)-\beta(\sigma\setminus\lambda)$.
\end{proof}

The final lemma of this section allows us to compare $\cha_\lambda(\mu)$ with $\cha_\lambda(g_k(\mu))$.

\begin{lemma}\label{signs2}
If $\lambda=(\lambda_0,...,\lambda_{n-1})$ with $|\lambda_j|=d_j$, then
\[
\cha_\lambda(g_k(\mu))=\xi_n^{-k\sum jd_j}\cha_\lambda(-\mu).
\]
\end{lemma}

\begin{proof}
Write $\mu=(\mu^0,...,\mu^{n-1})$ with $\mu^s=(d_1^s,...,d_{l_s}^s)$ and define $(,)$ to be the inner product for which $\{v_\lambda\}$ is an orthonormal basis.  By \eqref{eqn:create2}, we have
\begin{align*}
\cha_\lambda(g_k(\mu))&=\left(\prod_{s=0}^{n-1}\prod_{i=0}^{l_s}\left( \sum_{j=0}^{n-1}\xi_n^{-d_i^skj+sj}\alpha_{-d_i^s}^j \right)v_\emptyset,v_\lambda\right)\\
&=\xi_n^{-k\sum jd_j}\left(\prod_{s=0}^{n-1}\prod_{i=0}^{l_s}\left( \sum_{j=0}^{n-1}\xi_n^{sj}\alpha_{-d_i^s}^j \right)v_\emptyset,v_\lambda\right)\\
&=\xi_n^{-k\sum jd_j}\cha_\lambda(-\mu).
\end{align*}

\end{proof}

\section{Proof of Theorem \ref{thm1}}\label{sec:proof}

We now check identities (I) - (III) of Reduction \ref{reduction}.

\subsection*{Identity (I)}

This follows immediately from Lemma \ref{lem:frame}.

\subsection*{Identity (II)}

Since $z_\mu=z_{\mu^0}z_{\mu^{tw}}$, we must show that
\[
\sum_{|\lambda|=|\mu|} \tilde{P}_\lambda(0)\cha_\lambda(\mu)=\left(\prod_{j=1}^{l_0}\frac{\sqrt{-1}(-1)^{d_j^0}}{2}\csc\left( \frac{d_{j}^0u}{2}\right)\right)\left( \sum_{|\sigma|=|\mu^{tw}|} \tilde{P}_{\sigma}(0)\cha_{\sigma}(\mu^{tw}) \right).
\]
after the change of variables.  To do this, it is equivalent to show
\[
\sum_{|\lambda|=|\mu|+k} \tilde{P}_\lambda(0)\cha_\lambda(\mu\cup \{k\}))=\frac{\sqrt{-1}(-1)^k}{2}\csc\left( \frac{ku}{2}\right)\left( \sum_{|\sigma|=|\mu|} \tilde{P}_{\sigma}(0)\cha_{\sigma}(\mu) \right)
\]
which is equivalent (before the change of variables) to
\begin{equation}\label{identii}
\sum_{|\lambda|=|\mu|+k} \tilde{P}_\lambda(0)\cha_\lambda(\mu\cup \{k\})= \frac{(-1)^kq^{\frac{k}{2}}}{1-q^k}\sum_{|\sigma|=|\mu|} \tilde{P}_{\sigma}(0)\cha_{\sigma}(\mu).
\end{equation}

Fix $\sigma$.  Then
\begin{align}\label{blah}
\nonumber\frac{(-1)^kq^{\frac{k}{2}}}{1-q^k}\tilde{P}_{\sigma}(0)\cha_{\sigma}(\mu)&=\frac{(-1)^{k+|\mu|}q^{\frac{k}{2}}}{1-q^k}\cha_{\sigma}(\mu)\frac{\cha_{\bar\sigma}(n^{|\mu|})}{\dim(\sigma)}q^\frac{|\mu|}{2}S_{\sigma}\\
&\nonumber=(-1)^{k+|\mu|}q^\frac{|\mu|+k}{2}\cha_{\sigma}(\mu)\frac{\cha_{\bar\sigma}(n^{|\mu|})}{\dim(\sigma)}\sum_{\bar\lambda\supset\bar\sigma}(-1)^{ht(\bar\lambda\setminus\bar\sigma)-1}S_\lambda\\
&\nonumber=\cha_{\sigma}(\mu)\sum_{\bar\lambda\supset\bar\sigma}(-1)^{\beta(\lambda\setminus\sigma)}\frac{\cha_{\bar\lambda}(n^{|\mu|+k})}{\dim(\lambda)}q^\frac{|\lambda|}{2}(-1)^{|\lambda|}S_\lambda\\
&=\cha_{\sigma}(\mu)\sum_{\bar\lambda\supset\bar\sigma}(-1)^{\beta(\lambda\setminus\sigma)}\tilde{P}_\lambda(0).
\end{align}
where the sum is over all $\bar\lambda$ obtained from $\bar\sigma$ by adding a length $kn$ border strip.  The first equality follows from Corollary \ref{sp}, the second from Theorem \ref{thm:comb1}, the third from Lemma \ref{signs1}, and the fourth is another application of Corollary \ref{sp}.

From (\ref{eqn:create2}), we know
\begin{equation}\label{signs3}
\cha_\lambda(\mu\cup \{k\})=\sum_{\sigma}\cha_\sigma(\mu)(-1)^{\beta(\lambda\setminus\sigma)},
\end{equation}
where the sum is over all $\sigma$ such that $\bar\sigma$ is obtained from $\bar\lambda$ by removing a $kn$ strip. Summing \eqref{blah} over all $\sigma$ proves identity \eqref{identii} and thus (II).

\subsection*{Identity (III)}

Applying Lemma \ref{signs2}, (III) is equivalent to
\[
\sum_\nu\left( \sum_\lambda \tilde{P}_\lambda(0)\frac{\cha_\lambda(\nu)}{z_\nu} \right)z_\nu \left( \sum_\sigma \xi_n^{-k\sum j|\sigma_j|}\frac{\cha_\sigma(-\nu)}{z_\nu}\frac{\cha_\sigma(\mu)}{z_\mu}\re^{\frac{k}{n}(\sqrt{-1}f_T(\sigma)u+\sum\xi_{2n}^{-i} f_{i}(\sigma)x_i)}\right)=0.
\]
Summing over all $\nu$ and using orthogonality of characters, the left side becomes
\[
\sum_\lambda \tilde{P}_\lambda(0)\frac{\cha_\lambda(\mu)}{z_\mu}\xi_n^{-k\sum j|\lambda_j|}\re^{\frac{k}{n}(\sqrt{-1}f_T(\lambda)u+\sum \xi_{2n}^{-i}f_{i}(\lambda)x_i)}.
\]
Applying Lemma \ref{lem:frame}, we then see that (III) is equivalent to
\begin{equation*}
\sum_\lambda \tilde{P}_\lambda(0)\cha_\lambda(\mu)\left(\prod_{(i,j)\in\bar\lambda} q_{j-i}^{j-i}\right)^{k/n}\hspace{-.5cm}=0
\end{equation*}
for any $\mu$ with at least one untwisted part.  This is equivalent to
\begin{equation}\label{identiii}
\sum_\lambda \tilde{P}_\lambda(0)\cha_\lambda(\nu\cup\{k\})\left(\prod_{(i,j)\in\bar\lambda} q_{j-i}^{j-i}\right)^{k/n}\hspace{-.5cm}=0
\end{equation}
for any $\nu$.  Fix $\sigma$ with $|\sigma|=|\nu|$.  Then
\begin{align*}
0&=\sum_{\bar\lambda\supset\bar\sigma} (-1)^{ht(\bar\lambda\setminus\bar\sigma)-1}S_\lambda^k\\
&=\sum_{\bar\lambda\supset\bar\sigma}\frac{\cha_{\bar\sigma}(n^{|\sigma|})}{\dim(\sigma)}\cha_\sigma(\nu) (-1)^{ht(\bar\lambda\setminus\bar\sigma)-1}S_\lambda^k\\
&=\sum_{\bar\lambda\supset\bar\sigma}(-1)^{\beta(\lambda\setminus\sigma)}\cha_\sigma(\nu)q^{-\frac{|\lambda|}{2}}(-1)^{|\lambda|}\tilde{P}_\lambda(0)\left(\prod_{(i,j)\in\bar\lambda} q_{j-i}^{j-i}\right)^{k/n}.
\end{align*}
where the first equality is Theorem \ref{thm:comb2}, the second holds because $\sigma$ is fixed, and the third follows from Lemmas \ref{sp2} and \ref{signs1}.  Since $|\lambda|$ is constant over the sum, it follows that
\[
0=\cha_\sigma(\nu)\sum_{\bar\lambda\supset\bar\sigma}(-1)^{\beta(\lambda\setminus\sigma)}\tilde{P}_\lambda(0)\left(\prod_{(i,j)\in\bar\lambda} q_{j-i}^{j-i}\right)^{k/n}.
\]
Summing over all $\sigma$ (using equation \eqref{signs3}) proves \eqref{identiii} and thus finishes the proof of Theorem \ref{thm1}.

\section{GW/DT for local $\Z_n$-gerbes over $\proj^1$}\label{sec:app}

We conclude by giving an application of the gerby Gopakumar-Mari\~no-Vafa formula.  In particular, we prove that the Gromov-Witten potential of any local $\Z_n$-gerbe over $\proj^1$ is equal to the reduced, multi-regular Donaldson-Thomas potential after an explicit change of variables.

\begin{definition}
A \textit{local $\Z_n$-gerbe over $\proj^1$} is the total space of a rank two Calabi-Yau bundle $L_1\oplus L_2$ over some $\G_k$ with trivial generic isotropy.
\end{definition}

The CY condition implies that $\deg(L_1)+\deg(L_2)=-2$.  Because of the generically trivial isotropy, we know that the $\Z_n$ isotropy acts on the fibers of $L_1$ by a generator $\zeta\in{\Z_n}$ and on the fibers of $L_2$ by its inverse $\zeta^{-1}$.  The automorphism of $\Z_n$ which maps $\zeta\rightarrow\xi$ induces an isomorphism of the total space which allows us to assume that the isotropy always acts on the fibers of $L_1$ and $L_2$ with weights $\xi$ and $\xi^{-1}$, respectively (c.f. discussion after Definition \ref{gerbe}).

Fix $k\in\{0,...,n-1\}$ and set $e:=\gcd(k,n)$.  Then $\Pic(\G_k)=\frac{e}{n}\Z\oplus\Z_e$.  For each $b\in\frac{e}{n}\Z\oplus\Z_e$ we let $\La_b$ denote the corresponding orbifold line bundle.  The subset of $\Pic(\G_k)$ where $\Z_n$ acts on fibers as multiplication by $\xi$ is given by $(\Z-\frac{k}{n})\oplus\{1\}$.  Every local $\Z_n$-gerbe over $\proj^1$ is isomorphic to $\X_{k,b}:=\Tot\left(\La_b\oplus\La_{-b-2}\right)$ for some $k\in\{0,...,n-1\}$ and $b\in\Z-\frac{k}{n}$.

By the gluing formula of \cite{r:lgoa}, the degree $d$ Gromov-Witten potential of $\X_{k,b}$ is given by
\begin{equation}\label{gwdt1}
\hspace{-.75cm}GW_d(\X_{k,b})=\sum_\mu V_\mu^\bullet\left(b\right) z_\mu V_{g_k(\mu)}^\bullet\left(0\right) \prod_{i,j}(-1)^{d_j^ib+1+\delta_{0,i}+\delta_{0,(-d_j^ik-i) \text{mod} n}+\frac{i}{n}+\frac{(d_j^ik-i) \text{mod} n}{n}}.
\end{equation}
where the sign is the gluing term in \cite{r:lgoa}.

Analyzing the modification in \eqref{modgw}, we see that \eqref{gwdt1} is equivalent to
\begin{equation}\label{gwdt2}
GW_d(\X_{k,b})=(-1)^{db}\sum_\mu \tilde{V}_\mu^\bullet\left(b\right) z_\mu \tilde{V}_{g_k(\mu)}^\bullet\left(0\right).
\end{equation}

Applying the change of variables in Theorem \ref{thm1}, then using Lemma \ref{signs2} and orthogonality of characters, we find that
\begin{align*}
GW_d(\X_{k,b})&=(-1)^{db}\sum_\mu\left( \sum_\lambda \tilde{P}_\lambda\left(b\right)\frac{\cha_\lambda(\mu)}{z_\mu} \right)z_\mu\left(\sum_\sigma \tilde{P}_\sigma\left(0 \right)\frac{\cha_\sigma(g_k(\mu))}{z_{g_k(\mu)}}\right)\\
&=(-1)^{db}\sum_\lambda\xi_n^{-k\sum i|\lambda_i|}\tilde{P}_\lambda\left( b \right)\tilde{P}_\lambda\left( 0 \right)
\end{align*}

From equation \eqref{moddt}, we see that this last expression is
\begin{equation}\label{gwdt3}
\sum_\lambda P_{\lambda}(q_0,q_1,...,q_{n-1})E_{\lambda}P_{\lambda'}(q_0,q_{n-1},...,q_1)
\end{equation}
where
\[
E_{\lambda}:=\prod_{(i,j)\in\bar\lambda}q_{j-i}^{(b+2)i-bj-1}(-1)^{dnb}.
\]
By the main result of \cite{bcy:otv}, \eqref{gwdt3} is equal to the reduced, multi-regular, degree $d$ Donaldson-Thomas potential $DT_{mr,d}'(\X_{k,b})$ after the substitution $q_0\rightarrow-q_0$.  This proves Theorem \ref{thm2}.

\bibliographystyle{gtart}

\end{document}